\documentclass[11pt, reqno]{amsart}
\usepackage{lmodern}
\usepackage{amsmath, amsthm, amssymb, amsfonts}
\usepackage[normalem]{ulem}
\usepackage{hyperref}
\usepackage[all,cmtip]{xy}
\usepackage{verbatim}
\usepackage{nccmath}
\usepackage{stmaryrd}
\usepackage{caption}
\usepackage{mathrsfs}
\usepackage{mathtools}
\usepackage{esvect}
\usepackage{cite}
\usepackage{bbm}
\usepackage{eucal}

\usepackage{mathbbol}
\usepackage{tabularx}
\usepackage[toc,page]{appendix}
\usepackage{tikz-cd}

\newcommand{\BC}{{\mathbb{C}}}

\newcommand{\BE}{{\mathbb{E}}}
\newcommand{\BF}{{\mathbb{F}}}

\newcommand{\BH}{{\mathbb{H}}}

\newcommand{\BL}{{\mathbb{L}}}

\newcommand{\BQ}{{\mathbb{Q}}}
\newcommand{\BR}{{\mathbb{R}}}

\newcommand{\BT}{{\mathbb{T}}}

\newcommand{\BZ}{{\mathbb{Z}}}

\newcommand{\CC}{{\mathcal C}}

\newcommand{\CF}{{\mathcal F}}

\newcommand{\CH}{{\mathcal H}}

\newcommand{\CL}{{\mathcal L}}
\newcommand{\CM}{{\mathcal M}}
\newcommand{\CN}{{\mathcal N}}
\newcommand{\CO}{{\mathcal O}}
\newcommand{\CP}{{\mathcal P}}

\newcommand{\CS}{{\mathcal S}}

\newcommand{\CY}{{\mathcal Y}}

\newcommand{\Mbar}{{\overline M}}
\newcommand{\p}{\mathbb{P}}
\newcommand\Spec{\operatorname{Spec}}
\newcommand\Ext{\operatorname{Ext}}
\newcommand{\tr}{{\mathrm{tr}}}
\DeclareMathOperator{\Hilb}{Hilb}
\newcommand{\rk}{{\mathrm{rk}}}
\newcommand{\ch}{{\mathrm{ch}}}
\newcommand{\rCoh}{\mathfrak{Coh}_{r}}
\newcommand{\ev}{{\mathrm{ev}}}
\newcommand{\CHom}{{\mathcal{H} om}}
\newcommand{\Pic}{\mathop{\rm Pic}\nolimits}
\newcommand{\rCohS}{\mathfrak{Coh}_{r}(S)}
\newcommand{\Cohc}{\mathrm{Coh}}

\newcommand{\FM}{\mathfrak{M}}
\newcommand{\Eff}{\mathrm{Eff}}

\newcommand{\rpCoh}{\mathfrak{Coh}_{r}^{\sharp}}
\newcommand{\Coh}{\mathfrak{Coh}}

\newtheorem{theorem}{Theorem}[section]
\newtheorem{lemma}[theorem]{Lemma}
\newtheorem{proposition}[theorem]{Proposition}
\newtheorem{corollary}[theorem]{Corollary}
\theoremstyle{definition}
\newtheorem{definition}[theorem]{Definition}
\newtheorem{remark}[theorem]{Remark}

\numberwithin{equation}{section}
\newtheorem*{claim}{Claim}

\begin{document}

		\title[Quasimaps to moduli spaces of  sheaves on a $K3$ surface]{Quasimaps to moduli spaces of  sheaves on a $K3$ surface}

			\author{Denis Nesterov}
		\address{ETH Z\"urich, Departement Mathematik}
		\email{denis.nesterov@math.ethz.ch}
		\maketitle

		\begin{abstract} In this article, we study quasimaps to moduli spaces of sheaves on a $K3$ surface $S$.  We construct a surjective cosection of the obstruction theory of moduli spaces of quasimaps. We then establish reduced wall-crossing formulas which relate the reduced Gromov--Witten theory of moduli spaces of sheaves on $S$ and the reduced Donaldson--Thomas theory of $S\times C$, where $C$ is a nodal curve. As applications, we prove the Hilbert-schemes part of the Igusa cusp form conjecture; higher-rank/rank-one  Donaldson--Thomas correspondence with relative insertions on $S\times C$, if $g(C)\leq1$;  Donaldson--Thomas/Pandharipande--Thomas correspondence with relative insertions on $S\times \p^1$.
		\end{abstract}

 \setcounter{tocdepth}{1}
		\tableofcontents

\section{Introduction}
	
	\subsection{Overview} In \cite{N},  $\epsilon$-stable quasimaps to a moduli space of sheaves on a surface $S$ were introduced. When applied to a  Hilbert scheme of points $S^{[n]}$, moduli spaces of $\epsilon$-stable quasimaps $Q^{\epsilon}_{g,N}(S^{[n]},\beta)$ interpolate between moduli spaces of stable maps to $S^{[n]}$ and moduli spaces of 1-dimensional subschemes on $S\times C$ for a moving nodal curve $C$, 
	\begin{equation} \label{qmwall1}
		\xymatrix{
			\Mbar_{g,N}(S^{[n]},\beta)  \ar@{<-->}[r]|-{\epsilon } &\mathrm{Hilb}_{n,\check{\beta}}(S\times C_{g,N})}.
	\end{equation}
	Zhou's master-space technique \cite{YZ} leads to quasimap wall-crossing formulas associated to a variation of the stability parameter $\epsilon\in \BR_{>0}$.  By (\ref{qmwall1}), these wall-crossing formulas therefore relate the Gromov--Witten (GW) theory of $S^{[n]}$  and the Donaldson--Thomas (DT) theory of $S\times C$ with relative insertions.

	The case of moduli spaces of sheaves on a $K3$ surface requires a special treatment due to the presence of holomorphic symplectic forms and, consequently, the  vanishing of the standard virtual fundamental class on 	$\Mbar_{g,N}(S^{[n]},\beta)$. In more precise terms, the vanishing is due to existence of a surjective cosection of the obstruction-theory complex $\BE^{\bullet}$,
	\[\sigma \colon \BE^{\bullet} \twoheadrightarrow \CO[-1].\]
	A non-trivial \textit{reduced} enumerative theory is obtained by taking the cone of  $\sigma$.  The same phenomenon happens on the DT side. The obstruction theory of moduli spaces  $\mathrm{Hilb}_{n,\check{\beta}}(S\times C_{g,N})$ admits a surjective cosection, hence the reduction is also necessary.  In order to relate the reduced GW theory of $S^{[n]}$ and the reduced DT theory of $S\times C$, we have to  furnish $Q^{\epsilon}_{g,N}(S^{[n]},\beta)$ with a surjective cosection, and, consequently, with a reduced obstruction theory.  This is the principle aim of the present work.
	
	Once the reduced wall-crossing formula is established in Theorem \ref{reducedws}, we proceed to proving the following results:
	\begin{itemize}
		\item  the (reduced) quantum cohomology of $S^{[n]}$ is determined by  the relative Pandharipande--Thomas (PT) theory of $S\times \p^1$, if $S$ is a $K3$ surface, conjectured in \cite{OP1};
		\item the Hilbert-schemes part of the Igusa cusp form conjecture, conjectured in \cite{OP10b};
		\item  higher-rank/rank-one DT correspondence with relative insertions for  $S \times \p^1$ and $S \times E$, if $S$ is a $K3$ surface and $E$ is an elliptic curve;
		\item  DT/PT  correspondence with relative insertions for  $S \times \p^1$, if $S$ is a $K3$ surface.
	\end{itemize}
	\subsection{Cosection}  
	Let $S$ be a $K3$ surface, and $M(\mathbf{v})$ be a projective moduli space of stable sheaves on $S$ in a class $\mathbf{v} \in K_{\mathrm{num}}(S)$. To give a short motivation for our forthcoming considerations, let us recall the origin of reduced perfect obstruction theory of the GW theory of $M(\mathbf{v})$. Since $M(\mathbf{v})$ is hyper-K\"ahler, for any algebraic curve class $\beta \in H_2(M(\mathbf{v}),\BZ)$, there exists a deformation of $M(\mathbf{v})$ over a small disk $\Delta$, 
	\[\CM \rightarrow\Delta,\] for which the horizontal lift of $\beta$ is of $(k,k)$ type only at the central fiber (for our purposes, $\Delta$ can be replaced by a first-order neighbourhood of the origin $\Spec (\BC[\epsilon]/\epsilon^2)$). In particular, the standard GW invariants vanish.  We will call such family a \textit{twistor family} associated to the class $\beta$, and refer to \cite[Section 2.1]{KTred} for more about twistor families in the context of enumerative geometry. 
	
	To get a non-trivial enumerative theory, we have to remove obstructions that arise via such deformations of $M(\mathbf{v})$. However, in the case of $\epsilon$-stable quasimaps, we need twistor families not just  of the moduli space $M(\mathbf{v})$ but of the pair $(M(\mathbf{v}),\rCoh(S,\mathbf{v}))$, where $\rCoh(S,\mathbf{v})$ is the rigidified stack of all sheaves in the given class \cite{N}. Such twistor families can be given by non-commutative deformations of $S$. Let us now elaborate on this point by slightly changing the point of view.

	For simplicity, assume $M(\mathbf{v})=S^{[1]}=S$. A map $f\colon C \rightarrow S$ of degree $\beta$ is determined by its graph on $S\times C$. Let $I$ be the associated ideal sheaf of this graph. Deformation theories of $I$, as a sheaf with fixed determinant, and of $f$ are equivalent. Assuming $C$ is smooth and $\beta\neq 0$, the existence of a first-order twistor family associated to the class $\beta$ is therefore equivalent to the surjectivity of the following composition 
	\begin{equation} \label{semiregularity}
		H^1(T_S) \hookrightarrow H^1(T_{S\times C}) \xrightarrow{\cdot \mathrm{At}(I) } \Ext^2(I,I)_0 \xrightarrow{\sigma_I} H^3(\Omega_{S\times C}^1)\cong \BC,
	\end{equation}
	i.e.\ to the existence of a class $\kappa_{\beta} \in H^1(T_S)$ whose image is non-zero with respect to the composition above, where $\sigma_I:=\tr(* \cdot -\mathrm{At}(I))$ for the Atiyah class $\mathrm{At}(I) \in \Ext^1(I,I\otimes \Omega^1_{S\times C})$.  To see this, recall that the second map gives the obstruction to deform $I$ along a first-order deformation $\kappa \in H^1(T_S)$, while the third map, called \textit{semiregularity map} \cite{BF}, relates obstructions of deforming $I$ to the obstructions of
	\[\ch_2(I)=(-\beta, -1) \in H^4(S\times C,\BZ)=H^2(S,\BZ)\oplus \BZ\]
	to stay of Hodge type $(k,k)$. With these interpretations in mind, it is not difficult to grasp that  $\kappa_{\beta}$ is indeed our first-order twistor family associated to $\beta$.

	The semiregularity map $\sigma_I$ globalises, i.e.\ there exists a cosection 
	\[\sigma \colon \BE^{\bullet} \twoheadrightarrow \CO[-1]\]of the obstruction-theory complex of the moduli space of ideals  $\Hilb(S\times C)$.   This cosection $\sigma$ is surjective by the existence of first-order twistor families, if the second Chern character of ideals is equal to $(\beta, n)$ for $\beta \neq 0$. By the cosection localisation introduced by  Kiem--Li  \cite{KiL}, the standard virtual fundamental class therefore vanishes. To make the enumerative theory non-trivial, we have to consider the reduced obstruction-theory complex $\BE^{\bullet}_{\mathrm{red}}:=\mathrm{cone}(\sigma)[-1]$. Proving that $\BE^\bullet_{\mathrm{red}}$ defines an obstruction theory,
	\[\BT_{\Hilb} \rightarrow\BE^\bullet_{\mathrm{red}}, \]
	is difficult, we show  it under a certain assumption in Proposition \ref{reduction}. Instead, \cite{KiL} provides a construction of the reduced virtual fundamental class without an obstruction-theory morphism.
	
	Let us come back to the case of a general moduli space of sheaves $M(\mathbf{v})$. By construction of $M(\mathbf{v})$, the deformation theory of quasimaps to $M(\mathbf{v})$ is equivalent to the one of sheaves on threefolds of the type $S\times C$, see  \cite[Proposition 5.1]{N} for more details. The obstruction theory of higher-rank sheaves on $S\times C$ also admits a cosection given by the semiregularity map. We want to show it is surjective. However, already for $M(\mathbf{v})=S^{[n]}$ with $n>1$, there is a problem with the argument presented above. If the degree of $f \colon C \rightarrow S^{[n]}$ is equal to a multiple of the exceptional curve class\footnote{The curve class dual to a multiple of the exceptional divisor associated to the resolution of singularities $S^{[n]} \rightarrow S^{(d)}$.}, then (\ref{semiregularity}) is zero. Indeed, in this case, $\ch_{2}(I)=(0,n)$, and the composition (\ref{semiregularity}) is equal to the contraction $\langle-,\ch_{2}(I) \rangle$, which therefore pairs trivially with classes in $H^1(T_S)$. The geometric interpretation of this phenomenon is that the exceptional curve class of $S^{[n]}$ stays of Hodge type $(k,k)$ along  commutative deformations of $S$, because punctorial Hilbert schemes deform to punctorial Hilbert schemes under commutative deformations of $S$. To fix the argument, we have to consider classes not just in $H^1(T_S)$, but in a larger space \[H^0(\wedge^2T_S)\oplus H^1(T_S)\oplus H^2(\CO_S),\]
	i.e.\ we have to consider non-commutative first-order twistor families to prove the surjectivity of the semiregularity map. 
	
	\subsection{Strategy}	
	For surjectivity of the semiregularity map, we will largely follow \cite[Section 4]{BF} and \cite[Proposition 11]{MPT} with few extra layers of complications. Firstly, since our threefold $S\times C$ might be singular (because $C$ is nodal),  we have to consider Atiyah classes valued in $\Omega^1_{S} \boxplus \omega_C$,
	\[\mathrm{At}_{\omega}(F)\in \Ext^1(F,F\otimes (\Omega^1_S \boxplus \omega_{C})),\]
	instead of $\Omega^1_{S} \boxplus \BL_C=\Omega^1_{S} \boxplus \Omega^1_C$, as the latter does not behave well under degenerations. Chern characters of sheaves are then defined via the Atiyah class of the form as above. Secondly, after establishing the expected correspondence between degrees of quasimaps and Chern characters of sheaves,  we allow contractions with classes in 
	\[H^0(\wedge^2T_S)\oplus H^1(T_S)\oplus H^2(\CO_S),\]
	instead of just $H^1(T_S)$, unlike in \cite[Section 4]{BF}.  Proposition \ref{nonc} is a vast extension of \cite[Proposition 11]{MPT} and implies surjectivity of the global semiregularity map for higher-rank sheaves, Corollary \ref{cosection}.

	Having constructed a surjective cosection of the obstruction theory, ideally one would like to reduce the obstruction theory. However, due to the involvement of non-commutative geometry in our considerations, we can reduce the obstruction theory only under a certain assumption, which is  not unnatural, see Proposition \ref{reduction} for more details. For that reason, we do not use our reduced obstruction theory for the construction of the reduced virtual fundamental class.  We instead choose to work with reduced classes of \cite{KiL}.

	\subsection{Applications of the quasimap wall-crossing} \label{applications} 
	\subsubsection{Enumerative geometry of $S^{[n]}$}
	Let $E$ be a fixed elliptic curve. A moduli space of non-constant maps from $E$ to $S^{[n]}$ up to translations by $E$ is of reduced virtual dimension 0. Hence the associated GW invariants do not require insertions, and we can define a virtual count of $E$ inside $S^{[n]}$. Applying the reduced wall-crossing to these GW invariants, we obtain that, up to the wall-crossing terms, they are equal to rank-one DT invariants on the Calabi--Yau threefold $S\times E$,  
	\begin{equation} \label{igus}
		\mathsf{GW}_{E}(S^{[n]})=\mathsf{DT_{\rk=1}}(S\times E)+\mathsf{Wall}.
	\end{equation}
	In \cite{Ob}, the wall-crossing terms are shown to be virtual Euler numbers of certain Quot schemes, which are computed for $S^{[n]}$, if $S$ is a $K3$ surface. Both theories in (\ref{igus})  are subjects of the
	\textit{Igusa cusp form conjecture} \cite[Conjecture A]{OPa}, which consists of two parts:
	\begin{itemize}
		\item expressing DT invariants of $S\times E$ in terms of the Igusa cusp form;
		\item expressing the difference of DT invariants of $S\times E$ and GW invariants of $S^{[n]}$ associated to $E$ in terms of an explicit correction term. 
	\end{itemize}
	The first part was proved in a series of papers \cite{OS2},  \cite{OPi}.  While the second part is therefore  a consequence of the reduced quasimap wall-crossing (\ref{igus}) together with the computations of \cite{Ob}. This completes the proof of the Igusa cusp form conjecture and provides an expression for GW invariants of $S^{[n]}$ associated to $E$. We refer to Section \ref{Sectiongenus1} for more details  about this conjecture.

	
	In a similar vein, in  \cite{Ob22}, a holomorphic anomaly equation is established for $S^{[n]}$ for genus-0 GW invariants with at most 3 markings. The proof crucially uses the quisimap wall-crossing, which relates genus-0 GW invariants of $S^{[n]}$ to PT invariants of $S\times \p^1$ and then to GW invariants of $S\times \p^1$ by \cite{Ob21}. 
	\subsubsection{Higher-rank/rank-one DT correspondence for $K3\times C$} Assume $M(\mathbf{v})$ satisfies various assumptions of  \cite{N} which are listed in Section \ref{surjcos}. In particular,  $M(\mathbf{v})$ is a smooth projective hyper-K\"ahler variety. Hence by \cite[Theorem 4.6]{HuyK}, it is deformation equivalent to $S^{[n]}$.  Using the quasimap wall-crossing on both sides, we can relate higher-rank DT theory to rank-one DT theory.  The pictorial representation of this procedure is given in Figure \ref{rk12}.

	If $(g,N)=(0,3)$, the rank-one DT side of the square is concerned with moduli spaces of 1-dimensional subschemes on $S\times \p^1$, which are transversal to the vertical divisors over $\{0,1,\infty\}\subset \p^1$, 
	\[S_{p}= S\times \{p\}\subset S\times \p^1, \quad p \in \{0,1,\infty\}.\]
	We allow $\p^1$ to sprout rational tails at $\{0,1,\infty\}$ in order to make the space proper.  The resulting space is denoted by  
	\begin{equation} \label{not}
		\mathrm{Hilb}_{n,\check{\beta}}(S\times \p^1/ S_{0,1,\infty}).
	\end{equation}
	Intersecting a vertical divisor with a subscheme, we obtain evaluation maps to $S^{[n]}$, 
	\[\ev_ p \colon \mathrm{Hilb}_{n,\check{\beta}}(S\times \p^1/ S_{0,1,\infty}) \rightarrow S^{[n]}, \quad  p \in \{0,1,\infty\}. \]
	Relative insertions are defined to be pullbacks of classes from $S^{[n]}$ via these evaluation maps. These are the insertions that can be compared to primary GW insertions in Figure \ref{rk12}. For the higher-rank DT side, there exist similar spaces, for the definition of which we refer to Section \ref{p1}.

	If $(g,N)=(0,3)$, the wall-crossing is trivial, hence Figure \ref{rk12}  gives us the following relation between DT invariants with relative insertions,
	\[\mathsf{DT_{rel,\rk=1}}(S\times \p^1/S_{0,1,\infty})=\mathsf{GW}_{0,3}(S^{[n]})=\mathsf{DT_{rel,\rk>1}}(S\times \p^1/S_{0,1,\infty}).\]
	In the case of $S\times E$, where $E$ is an elliptic curve, we also get a wall-crossing statement for absolute invariants and equality of certain relative invariants, as is explained in Section \ref{ell}.

	\vspace{1 in}
	\begin{figure}[h!] 
		\scriptsize
		\[
		\begin{picture}(200,75)(-30,-50)
			\thicklines
			\put(25,-25){\line(1,0){20}}
			\put(125,-25){\line(-1,0){20}}
			\put(25,-25){\line(0,1){40}}
			\put(25,30){\makebox(0,0){\textsf{quasimap}}}
			\put(25,20){\makebox(0,0){\textsf{wall-crossing}}}
			\put(25,35){\line(0,1){40}}
			\put(125,-25){\line(0,1){40}}
			\put(125,30){\makebox(0,0){\textsf{quasimap}}}
			\put(125,20){\makebox(0,0){\textsf{wall-crossing}}}
			\put(25,75){\line(1,0){25}}
			\put(100,75){\line(1,0){25}}
			\put(125,35){\line(0,1){40}}
			\put(140,85){\makebox(0,0){$\mathsf{GW}(S^{[n]})$}}
			\put(15,85){\makebox(0,0){$\mathsf{GW}(M(\mathbf{v}))$}}
			\put(75,80){\makebox(0,0){\textsf{deformation}}}
			\put(75,70){\makebox(0,0){\textsf{invariance}}}
			\put(5,-35){\makebox(0,0){$\mathsf{DT_{rel,\rk>1}}(S\times C_{g,N})$}}
			\put(150,-35){\makebox(0,0){$\mathsf{DT_{rel,\rk=1}}(S\times C_{g,N})$}}
			\put(75,-25){\makebox(0,0){$\mathsf{DT_{\rk>1}/DT_{\rk=1}}$}}
		\end{picture}
		\]
		\caption{Higher-rank/rank-one DT correspondence}
		\label{rk12}
	\end{figure}
	

	\subsubsection{\textsf{DT/PT} correspondence for $K3 \times C$} 
	By results from  \cite[Section 6.3]{N}, there also exists a theory of quasimaps that interpolates between moduli spaces of stable maps to $S^{[n]}$ and stable pairs on $S\times C$ in the sense of \cite{PT},  
	\begin{equation} \label{qmwall}
		\xymatrix{
			\Mbar_{g,N}(S^{[n]},\beta)  \ar@{<-->}[r]|-{\epsilon } &\mathrm{P}_{n,\check{\beta}}(S\times C_{g,N})}.
	\end{equation}
	The resulting wall-crossing formulas therefore relate the GW theory of $S^{[n]}$ to the PT theory of $S\times C$ with relative insertions. We refer to this wall-crossing as \textit{perverse} quasimap wall-crossing, for more details on this terminology we refer to \cite[ Section 6.3]{N}.

	Using both standard and perverse quasimap wall-crossings, we can reduce the  DT/PT  correspondence   for a relative geometry of the form 
	\[ S \times C_{g,N} \rightarrow \Mbar_{g,N}\]
	to  the DT/PT  correspondence of wall-crossing invariants, as is illustrated in Figure \ref{DT/PTdiagram}.	As before, if $(g,N)=(0,3)$,  the wall-crossing is trivial. We therefore obtain the following relation,
	\[\mathsf{DT_{rel,\rk=1}}(S\times \p^1/S_{0,1,\infty})=\mathsf{GW}_{0,3}(S^{[n]})=\mathsf{PT_{rel,\rk=1}}(S\times \p^1/S_{0,1,\infty}).\]
	Note that this is an equality of invariants, unlike the more conventional DT/PT correspondence for Calabi--Yau threefolds which involves dividing by $0$-dimensional DT invariants, as it was conjectured in \cite[Conjecture 3.3]{PT} and proved in \cite[Theorem 1.1]{Bri}.  This form of DT/PT correspondence is not surprising due to the nature of  \textit{reduced} virtual fundamental classes. Moreover, since  we are in the setting of a non-Calabi-Yau relative geometry,  the techniques of wall-crossings in derived categories of \cite{KS} and \cite{JS} cannot be applied to prove wall-crossing statements as above.
	
	\vspace{0.1in}
	\begin{figure}[h!] \vspace{0.1in}
		\scriptsize
		\[
		\begin{picture}(200,75)(-30,-50)
			\thicklines
			\put(60,10){\line(1,1){15}}
			\put(25,-25){\line(1,1){15}}
			\put(25,-25){\line(1,0){30}}
			\put(125,-25){\line(-1,0){30}}
			\put(125,-25){\line(-1,1){15}}
			\put(90,10){\line(-1,1){15}}
			\put(75,30){\makebox(0,0){$\mathsf{GW}(S^{[n]})$}}
			\put(50,5){\makebox(0,0){\textsf{quasimap}}}
			\put(50,-5){\makebox(0,0){\textsf{wall-crossing}}}
			\put(105,5){\makebox(0,0){\textsf{quasimap}}}
			\put(105,-5){\makebox(0,0){\textsf{wall-crossing}}}
			\put(5,-35){\makebox(0,0){$\mathsf{DT_{rel,\rk=1}}(S\times C_{g,N})$}}
			\put(150,-35){\makebox(0,0){$\mathsf{PT_{rel,\rk=1}}(S\times C_{g,N})$}}
			\put(75,-25){\makebox(0,0){$\mathsf{DT/PT}$}}
		\end{picture}
		\]
		\caption{DT/PT correspondence }
		\label{DT/PTdiagram}
	\end{figure}
\vspace{-0.5cm}
	\subsection{Notation and conventions}
	We work over the field of complex numbers $\BC$. We set 
	$e_{\BC^*}(\BC_{\mathrm{std}})=z,$
	where $\BC_{\mathrm{std}}$ is the weight 1 representation of $\BC^*$ on the vector space $\BC$. All functors are derived, unless stated otherwise. Cohomologies and homologies have rational coefficients, unless stated otherwise.

	\section{Semiregularity map}
	
	\subsection{Preliminaries}Let $S$ to be a $K3$ surface over the field of complex numbers $\BC$, and let $F$ be a sheaf on $S\times C$ flat over a nodal curve $C$. More generally, the following discussion also applies to perfect complexes with the help of \cite{HT}.  In particular, we may assume that $F$ is a stable pair in the sense of \cite{PT}, which is flat over nodes of $C$. 
	
	Consider the Atiyah class 
	\[\mathrm{At}(F) \in \Ext^1(F,F\otimes \Omega_{S\times C}^1),\] represented by the canonical exact sequence
	\[0\rightarrow F\otimes \Omega_{S\times C}^1 \rightarrow \CP^1(F) \rightarrow F \rightarrow 0,\]
	where $\CP^1(F)$ is the sheaf of principle parts. 
	Composing the Atiyah class with the natural map \[\Omega_{S\times C}^1=\Omega_S^1\boxplus \Omega_C^1 \rightarrow \Omega^1_S \boxplus \omega_{C},\]
	we obtain a class \[\mathrm{At}_{\omega}(F)\in \Ext^1(F,F\otimes (\Omega^1_S \boxplus \omega_{C})).\] We then define the Chern character of a sheaf $F$ on $S \times C$ for a possibly singular $C$ as follows,
	\begin{equation} \label{Chern}
		\ch_{k}(F):= \tr \left( \frac{(-1)^k}{k!}\mathrm{At}_{\omega}(F)^k \right) \in H^k(\wedge^k(\Omega^1_{S} \boxplus \omega_{C})).
	\end{equation}
	If $C$ is smooth, it agrees with the standard definition of the Chern character. Using the canonical identification $H^1(\omega_C)\cong \BC$, and
	\[\wedge^k(\Omega^1_{S}\boxplus \omega_C) \cong \Omega^{k}_{S}\boxplus (\Omega_S^{k-1}\boxtimes \omega_{S}),\]
	we get a K\"unneth's decomposition of the cohomology
	\[H^k(\wedge^k(\Omega^1_{S} \boxplus \omega_{C}))\cong H^k(\Omega^k_S)\oplus H^{k-1}(\Omega^{k-1}_S),\]
hence
	\begin{equation} \label{lambda}
		\bigoplus_k H^k(\wedge^k(\Omega^1_{S} \boxplus \omega_{C}))\cong \Lambda \oplus \Lambda(-1),
	\end{equation}
	where 	\begin{equation*}
		\Lambda:=\bigoplus_{k} H^{k,k}(S).
	\end{equation*}
	With respect to this decomposition above, the Chern character $\ch(F)$ has two components
	\[\ch(F)=(\ch(F)_{\mathrm{f}},\ch(F)_{\mathrm{d}})\in \Lambda \oplus \Lambda(-1),\]
	such that the first component is determined by the Chern character of a fiber of $F$ over a point $p\in C$, 
	\[ \ch(F)_{\mathrm{f}}=\ch(F_p).\]
	On the other hand, if $C$ is smooth, then $\ch(F)_{\mathrm{d}}$ is determined by the degree of the quasimap associated to $F$, see \cite[Lemma 3.3]{N} for more details. We will show that the definition of $\ch(F)$ in (\ref{Chern}) is compatible with the definition in \cite[Definition 3.5]{N} for a singular curve $C$. So let $C$ be singular, let 
	\[ \pi\colon  S\times \tilde{C} \rightarrow S \times C\] be
	the normalisation map, and $\pi^*F_{i}$ be the restriction of $\pi^*F$ to the connected components $\tilde{C}_{i}$ of $\tilde{C}$. The above decomposition of the Chern character  satisfies the following property. 
	\begin{lemma}\label{identification} Under the identification (\ref{lambda}), we have
		\[\ch(F)=(\ch(F_p),\sum_i \ch(\pi^*F_{i})_{\mathrm{d}}) \in \Lambda \oplus \Lambda(-1).\]
	\end{lemma}
	\textit{Proof.}
	Firstly, there exist canonical maps making the following diagram commutative,
	\[
	\begin{tikzcd}[row sep=small, column sep = small]
		0 \arrow[r]& \pi^*F\otimes \pi^*\Omega_{S\times C}^1 \arrow[r] \arrow[d] & \pi^*\CP^1(F) \arrow[r] \arrow[d] &\pi^*F \arrow[d,equal] \arrow[r] & 0 \\
		0 \arrow[r]& \pi^*F\otimes \Omega_{S\times \tilde{C}}^1 \arrow[r] & \CP^1(\pi^*F) \arrow[r] &\pi^*F \arrow[r]& 0
	\end{tikzcd}
	\] 
	where the first row is exact on the left, because $L\pi^*F\cong \pi^* F$, by  flatness\footnote{To see that, one can use the standard locally-free resolution of a flat sheaf; these resolutions are functorial with respect to pullbacks.} of $F$ over $C$. The diagram above implies that the pullback of the Atiyah class $\pi^*\mathrm{At}(F)$ is mapped to $\mathrm{At}(\pi^*F)$ with respect to the map
	\[\Ext^1(\pi^*F,\pi^*F \otimes \pi^*\Omega_{S\times C}^1) \rightarrow \Ext^1(\pi^*F,\pi^*F \otimes \Omega_{S\times \tilde{C}}^1).\]
	The same holds for $\pi^*\mathrm{At}^k(F)$. Consider now the following commutative diagram,
	\[
	\begin{tikzcd}[row sep=small, column sep = small]
		R\CHom(F,F \otimes \Omega_{S\times C}^{k}) \arrow[r] \arrow[d] & \Omega_{S\times C}^k  \arrow[dd] \arrow[rdd] & \\
		\pi_*R\CHom(\pi^*F, \pi^*F\otimes \pi^*\Omega_{S\times C}^k)  \arrow[d]  &   &\\
		\pi_*R\CHom(\pi^*F, \pi^*F\otimes \Omega_{S\times \tilde{C}}^k) \arrow[r]  & \pi_* \Omega_{S\times \tilde{C}}^k \arrow[r] & \wedge^k(\Omega^1_{S}\boxplus \omega_C) \\
	\end{tikzcd}
	\]
	such that the first vertical map is the composition
	\begin{multline*}
		R\CHom(F,F \otimes \Omega_{S\times C}^{k})\rightarrow \pi_*L\pi^* R\CHom(F,F\otimes \Omega_{S\times C}^{k})\\
		=\pi_*R\CHom(\pi^*F,\pi^*F\otimes L\pi^*\Omega_{S\times C}^k)
		\rightarrow \pi_* R\CHom(\pi^*F,\pi^*F\otimes \pi^*\Omega_{S\times C}^k),
	\end{multline*}
	where we used that $L\pi^*F\cong \pi^*F$. Taking the cohomology of the diagram above and using the exactness of $\pi_*$, we can therefore factor the map \[\Ext^k(F,F\otimes \Omega_{S\times C}^k) \rightarrow H^k(\wedge^k(\Omega^1_{S}\boxplus \omega_C))\]
	as follows,
	\begin{multline*}
		\Ext^k(F,F\otimes \Omega_{S\times C}^k) \rightarrow \Ext^k(\pi^*F,\pi^*F \otimes \pi^*\Omega_{S\times C}^k)\rightarrow \Ext^k(\pi^*F,\pi^*F \otimes \Omega_{S\times \tilde{C}}^k) \\ \rightarrow H^k( \Omega_{S\times \tilde{C}}^k )\cong H^k(\Omega_S^k)\oplus \bigoplus_i H^{k-1}(\Omega_S^{k-1})\otimes H^1(\omega_{\tilde{C}_i}) \\ \rightarrow H^k(\Omega_S^k)\oplus H^{k-1}(\Omega_S^{k-1})\otimes H^1(\omega_C)\cong H^k(\wedge^k(\Omega^1_{S}\boxplus \omega_C)).
	\end{multline*}
	With respect to the natural identifications $H^1(\omega_{\tilde{C}_i}) \cong \BC$ and $H^1(\omega_{C}) \cong \BC$, the last map in the sequence  becomes 
	\[ H^k(\Omega_S^k)\oplus \bigoplus_i H^{k-1}(\Omega_S^{k-1}) \xrightarrow{(\mathrm{id},+)} H^k(\Omega_S^k)\oplus H^{k-1}(\Omega_S^{k-1}).\]
	The claim then follows by tracking the powers of the Atiyah class $\mathrm{At}^k(F)$ along the maps above.   \qed 
	

\subsection{Semiregularity map}
By pulling back  classes in
\[HT^2(S):=H^0(\wedge^2T_S) \oplus  H^1(T_{S}) \oplus H^2(\CO_{S})\] 
from $S$ to $S\times C$, we will treat  $HT^2(S)$ as classes on $S\times C$. Let 
\[ \sigma_{i}:= \tr\left(* \cdot \frac{(-1)^i}{i!}\mathrm{At}_{\omega}(F)^i \right)\colon  \Ext^2(F,F) \rightarrow H^{i+2}(\wedge^i(\Omega^1_{S}\boxplus \omega_C))\] 
be a \textit{semiregularity map}.
\begin{lemma}\label{At}  The following diagram is commutative, 
\[
\begin{tikzcd}[row sep=small, column sep = small]
	H^{2-k}(\wedge^kT_S) \arrow[rr, "\cdot \frac{(-1)^k}{k!}\mathrm{At}_{\omega}(F)^k "] \arrow{dr}[swap]{\langle *, \ch_{k+i}(F) \rangle} && \Ext^2(F,F) \arrow[dl, "\sigma_i"] \\
	& H^{i+2}(\wedge^i(\Omega^1_{S}\boxplus \omega_C))
\end{tikzcd}
\]
\end{lemma}
\textit{Proof.}
If $i=0$, then $\sigma_{0}=\tr$ and the commutativity is implied by the following property of the contraction pairing,
\[\langle \kappa, \tr(\mathrm{At}_{\omega}(F)^k))\rangle =\tr\langle \kappa,\mathrm{At}_{\omega}(F)^k  \rangle. \]
The proof is presented in \cite[Proposition 4.2]{BF} for $k=1$, and is the same for other values of $k$.

If $i=1$, then for the commutativity of the digram we have to prove that 
\[\bigg\langle \kappa, \tr \left(\frac{\mathrm{At}_{\omega}(F)^{k+1}}{k+1!}\right) \bigg\rangle= \tr \left( \bigg\langle \kappa, \frac{\mathrm{At}^{k}_{\omega}(F)}{k!} \bigg\rangle\cdot \mathrm{At}_{\omega}(F) \right). \]
If $\kappa \in H^2(\CO_S)$, the equality follows trivially,  since there is no contraction. The case of $\kappa \in H^1(T_S)$ is treated in \cite[Proposition 4.2]{BF}. For $\kappa \in H^0(\wedge^2 T_S)$, we use the derivation property for contraction with a 2-vector field 
\[\langle \xi, \mathrm{At}^3_{\omega}(F) \rangle = 3 \langle \xi, \mathrm{At}^2_{\omega}(F) \rangle\cdot \mathrm{At}_{\omega}(F),\] 
which can be checked locally on a 2-vector field of the form $V\wedge W$. 
\qed
\\

Due to the decomposition
\[H^i(\wedge^i(\Omega^1_{S}\boxplus \omega_C))\cong H^i(\Omega^i_{S}) \oplus H^{i-1}(\Omega^{i-1}_S),\] 
there are two ways to contract a class in $H^i(\wedge^i(\Omega^1_{S}\boxplus \omega_C))$ with a class in $H^{2-k}(\wedge^kT_S)$: either via the first component of the decomposition above or via the second. Hence due to the wedge degree or the cohomological degree, only one component of $H^i(\wedge^i(\Omega^1_{S}\boxplus \omega_C))$ pairs non-trivially with $H^{2-k}(\wedge^kT_S)$ for a fixed $k$.   It is not difficult to check that contraction with the Chern character,
\[H^{2-k}(\wedge^kT_S) \xrightarrow{\langle-,\ch_{k+i}(F) \rangle} H^{i+2}(\wedge^i(\Omega^1_{S}\boxplus \omega_C)),\]
is therefore equal to $\langle-,\ch(F)_{\mathrm{f}} \rangle$ for $i=0$ and to $\langle-,\ch(F)_{\mathrm{d}} \rangle$ for $i=1$. Moreover, using the identification
\[H^{i+2}(\wedge^i(\Omega^1_{S}\boxplus \omega_C)) \cong H^2(\CO_S),\]
the contraction $\langle-,\ch(F)_{\mathrm{d/f}} \rangle$ with classes on $S\times C$ is identified with the contraction with classes on $S$.
\begin{proposition} \label{nonc} Assume \[\ch(F)_{\mathrm{f}} \wedge \ch(F)_{\mathrm{d}}\neq 0,\]
then there exists $\kappa \in HT^2(S)$, such that 
\[\langle \kappa, \ch(F)_{\mathrm{f}} \rangle=0 \quad \text{and} \quad 
\langle \kappa, \ch(F)_{\mathrm{d}} \rangle \neq 0.\]
In particular, the restriction of the semiregularity map to the traceless part of $\Ext^2(F,F)$,
\[\sigma_1\colon  \Ext^2(F,F)_{0} \rightarrow H^3(\Omega_S^1\boxplus \omega_C),\]
is non-zero.
\end{proposition}
\textit{Proof.} Using a symplectic form on $S$, we have the following identifications,
\[\wedge^2T_S\cong \CO_S, \quad T_S\cong \Omega_S^1, \quad \CO_{S}\cong \Omega_S^2.\] 
After applying these identifications and taking the cohomology, the pairing 
\begin{equation} \label{cont}
HT^2(S)\otimes H\Omega_{0}(S) \rightarrow H^2(\CO_S),
\end{equation}
which is given by the contraction of classes, becomes the intersection pairing
\[H\Omega_{0}(S)\otimes H\Omega_{0}(S) \rightarrow H^2(\Omega_S^2),\]
where $H\Omega_0(S)=\bigoplus_i H^i(\Omega^i)$. In particular, the pairing (\ref{cont}) is non-degenerate. We conclude that $\ch(F)_{\mathrm{d}}^{\perp}$ and $\ch(F)_{\mathrm{f}}^{\perp}$ are distinct, if and only if $\ch(F)_{\mathrm{d}}$ is not a multiple of $\ch(F)_{\mathrm{f}}$. Hence there exists a class $\kappa \in HT^2(S)$, such that 
\[\langle \kappa, \ch(F)_{\mathrm{f}} \rangle=0 \quad \text{and} \quad 
\langle \kappa, \ch(F)_{\mathrm{d}} \rangle \neq 0.\]
 By Lemma \ref{At} and the discussion afterwards, the property$\langle \kappa, \ch(F)_{\mathrm{f}} \rangle=0$ implies that 
\[ \kappa \cdot \exp(-\mathrm{At}_{\omega}(F) \in \Ext^2(F,F)_0,\] 
while the property $\langle \kappa, \ch(F)_{\mathrm{d}} \rangle=0$ implies that the restriction of the semiregularity map to $\Ext^2(F,F)_0$ is non-zero, as it is non-zero when applied to the element $\kappa \cdot \exp(-\mathrm{At}_{\omega}(F))$. 
\qed
\\

\begin{remark} From the point of view of quasimaps, the condition
\[\ch(F)_{\mathrm{f}} \wedge \ch(F)_{\mathrm{d}}\neq 0,\]
is equivalent to the fact that the quasimap $f\colon C \rightarrow \rCohS$ associated to $F$ is not constant. 
\end{remark}
The above result has a following geometric interpretation. With respect to the Hochschild--Kostant--Rosenberg (HKR) isomorphism, 
\[HT^2(S)\cong HH^2(S),\]
the space $HT^2(S)$ parametrises first-order non-commutative deformations of $S$, i.e.\ deformations of $\mathrm{D}^b(S)$. Given a first-order deformation $\kappa \in HT^2(S)$, the unique horizontal lift of $\ch(F)_{\mathrm{d/f}}$ relative to a Gauss--Manin connection associated to $\kappa$ should stay of Hodge type $(k,k)$, if and only if $\langle \kappa, \ch(F)_{\mathrm{d/f}} \rangle=0$. On the other hand, $\langle \kappa, \exp(-\mathrm{At}_{\omega}(F))$ gives an obstruction for deforming $F$ on $S\times C$ in the direction $\kappa$. Hence, by Lemma \ref{At}, the semiregularity map $\sigma_i$ relates obstructions to deform $F$ along $\kappa$ with obstructions of $\ch(F)_{\mathrm{d/f}}$ to stay of Hodge type $(k,k)$.  Proposition \ref{nonc} states that there exists a deformation $\kappa \in HT^2(S)$, for which $\ch(F)_{\mathrm{f}}$ stays of Hodge type $(k,k)$, but $\ch(F)_{\mathrm{d}}$ does not. From the point of view of quasimaps, this means that the moduli space of stable sheaves $M$ on $S$ associated to the class $\ch(F)_{\mathrm{f}}$  deforms along $\kappa$, but the quasimap associated to $F$ does not. 

For example, let $S$ be a $K3$ surface associated to a cubic fourfold $Y$, such that the Fano variety of lines $F(Y)$ is isomorphic to $S^{[2]}$. Then if we deform $Y$ away from the Hassett divisor  \cite{Ha}, $F(Y)$ deforms along, but the point class of $S$ does not. Therefore such deformation of $Y$ will give the first-order non-commutative deformation $\kappa \in HT^2(S)$ of $S$, such that $\ch(F)_{\mathrm{f}}=(1,0,-2)$ stays of Hodge type $(k,k)$, but $\ch(F)_{\mathrm{d}}=(0,0,k)$ does not. Note that  $\ch(F)_{\mathrm{d}}=(0,0,k)$ corresponds to multiplies of the exceptional curve class
in $S^{[2]}$. Indeed, there are no commutative deformations of $S$ that will make $(0,0,k)$ non-Hodge, because, in this case, the exceptional divisor deforms along with $S^{[2]}$. 

\section{Reduced wall-crossing} 
\subsection{Surjective cosection} \label{surjcos}
We fix  a very ample line bundle $\CO_{S}(1)\in \Pic(S)$, a class $\mathbf{v} \in K_{\mathrm{num}}(S)$, and another class $\mathbf{u} \in K_{0}(S)$, such that:
\begin{itemize}
\item $\rk(\mathbf{v})>0$,
\item  $\chi(\mathbf{v}\cdot \mathbf{u})=1$,
\item for $\mathbf{v}$ and $\CO_S(1)$,  all semistable\footnote{Semistablity is defined with respect to Hilbert polynomials.} sheaves are stable.
\end{itemize}
Let $M(\mathbf{v})$ be the moduli space of  $\CO_S(1)$-stable sheaves on $S$ in the class $\mathbf{v}$. The second assumption  implies that $M(\mathbf{v})$ is a fine moduli space, while the last assumption implies that $M(\mathbf{v})$ is smooth and projective. The first two assumptions are made for technical reasons, and, in principle, can be dropped. We refer to \cite[Section 1.6]{N} for a more detailed discussion about why these assumptions are made. 

\begin{remark} \label{example} An example of a moduli space $M(\mathbf{v})$ which satisfies the assumptions above will be a moduli space of sheaves in the class $\ch(\mathbf{v})=(2,\alpha,2k+1)$ for a polarisation such that $\deg(\alpha)$ is odd (or a generic polarisation that is close to a polarisation for which $\deg(\alpha)$ is odd). Firstly, $\rk(\mathbf{v})$ and $\deg(\mathbf{v})$ are coprime, therefore there are no strictly slope semistable sheaves. The class $\mathbf{u}=[\CO_S]-(k+2)[\CO_{\mathrm{pt}}]\in K_{0}(S)$ has the property $\chi(\mathbf{v}\cdot \mathbf{u})=1$. Moreover, \cite[Corollary A.7]{N} holds in this case, therefore the  space $M_{\mathbf{v},\check{\beta}}(S\times C)$ is a moduli space of stable sheaves for some \textit{suitable} polarisation. More specifically, such set-up can be arranged on an elliptic $K3$ surface.  
\end{remark} 

By \cite[Theorem 3.16]{N}, there exists an identification between a space of $\epsilon$-stable quasimaps  $Q^{\epsilon}_{g,N}(M(\mathbf{v}), \beta)$ and a certain relative moduli space of sheaves  $M^{ \epsilon}_{\mathbf{v},\check{\beta}}(S \times C_{g,N})$,
\begin{equation} \label{idents}
Q^{\epsilon}_{g,N}(M(\mathbf{v}), \beta) \cong M^{ \epsilon}_{\mathbf{v},\check{\beta}}(S \times C_{g,N}),
\end{equation}
such that the corresponding obstructions theories are isomorphic. The product $S\times C_{g,N}$ stands for the relative geometry 
\[S\times C_{g,N} \rightarrow \Mbar_{g,N},\] where $C_{g,N}$ is the universal curve over  $\Mbar_{g,N}$. This identification depends on the choice of the class $\mathbf{u} \in K_0(S)$, which we suppress from the notation. For $\epsilon=0^+$,  the $\BC$-valued points of $M^{0^+}_{\mathbf v, \check{\beta}}(S\times C_{g,N})$ are triples  
\begin{equation} \label{triples}(C, \mathbf{p}, F),
\end{equation}
such that: 
\begin{itemize}
\item $(C, \mathbf{p})$ is a prestable nodal curve (no rational tails),
\item a sheaf $F$ on $S\times C$ flat over $C$,
\item$ \ch(F)=(\ch(\mathbf v), \check \beta)$,
\item $F_p$ is stable for a general $p\in C$,
\item $F_p$ is stable, if $p$ is a node or a marking,
\item the group of automorphisms of $(C, \mathbf{p})$ fixing $F$ is finite,
\item $\det(p_{C*}(p_{S}^{*}\mathbf{u}\otimes F))\cong\CO_{C}$.
\end{itemize}
We will frequently use the identification (\ref{idents}) to transfer various constructions from  sheaves to quasimaps, and vice versa. If $M(\mathbf{v})=S^{[n]}$, then one can consider a moduli space of perverse quasimaps $Q^{\epsilon}_{g,N}(M(\mathbf{v}), \beta)^{\sharp}$, which admits an identification with a relative moduli space stable pairs of \cite{PT}, 
\[ Q^{\epsilon}_{g,N}(S^{[n]}, \beta)^{\sharp}\cong \mathrm{P}^{\epsilon}_{n,\check{\beta}}(S\times C_{g,N}).\] 
The following discussion applies to both kinds of quasimaps. 

Let 
\begin{align*}
&\pi\colon S\times \CC^{\epsilon}_{g,N}  \rightarrow M^{ \epsilon}_{\mathbf{v},\check{\beta}}(S \times C_{g,N}), \\
&\BF \in \Cohc(S\times \CC^\epsilon_{g,N})
\end{align*}
be the universal threefold and the universal sheaf of $M^{ \epsilon}_{\mathbf{v},\check{\beta}}(S \times C_{g,N})$, respectively.  Let 
\[ \tr\colon R\CH om_{\pi} (\BF, \BF) \rightarrow R\pi_*(\CO_{S\times \CC^\epsilon_{g,N}}) \]
be the universal trace map. The complex 
\[\BE^{\bullet}:=R\CH om_{\pi}(\BF, \BF)_0[1]=\mathrm{Cone}(\tr)\]
defines a perfect obstruction theory on $M^{ \epsilon}_{\mathbf{v},\check{\beta}}(S \times C_{g,N})$ relative to the moduli stack of nodal curves $\FM_{g,N}$. We construct a cosection of the obstruction theory via the global relative semiregularity map,
\[\mathsf{sr}\colon  \BE^{\bullet} \rightarrow R^3\pi_{*}(\Omega^1_{S}\boxplus \omega_{\CC^\epsilon_{g,N}})[-1],\]
and since
\[R^3\pi_{*}(\Omega^1_{S}\boxplus \omega_{\CC^\epsilon_{g,N}})\cong H^{2}(\CO_S) \otimes \CO_{M^{ \epsilon}_{\mathbf{v},\check{\beta}}(S \times C_{g,N})},\]
we have
\begin{equation} \label{relativecosection}
\mathsf{sr}\colon  \BE^{\bullet}\rightarrow H^{2}(\CO_S) \otimes \CO_{M^{ \epsilon}_{\mathbf{v},\check{\beta}}(S \times C_{g,N})}[-1]\cong  \CO_{M^{ \epsilon}_{\mathbf{v},\check{\beta}}(S \times C_{g,N})}[-1].
\end{equation}
By the identification (\ref{idents}), we get a cosection for the obstruction theory of  $Q^{\epsilon}_{g,N}(M(\mathbf{v}), \beta)$. Surjectivity of the cosection follows from the preceding results. 

\begin{corollary} \label{cosection} Assuming $\beta \neq 0$, the semiregularity map $\mathsf{sr}$ is surjective.
\end{corollary}
\textit{Proof.} Under the given assumption, the surjectivity of the semiregularity map $\mathsf{sr}$ follows from Proposition \ref{nonc} and Lemma \ref{identification}.
\qed
\\

Consider now the composition
\begin{equation} \label{comp}
\Ext_C^1(\Omega_C, \CO_C(-\mathbf{p})) \rightarrow \Ext^2(F,F)_0 \xrightarrow{\sigma_1} H^3(\Omega^1_S\boxplus \omega_C),
\end{equation}
where the first map defined by the following composition,
\[ \Ext_C^1(\Omega_C, \CO_C(-\mathbf{p})) \rightarrow \Ext_C^1(\omega_C, \CO_C) \xrightarrow{\cdot-\mathrm{At}_{\omega}(F)} \Ext^2(F,F)_0.\]
The composition (\ref{comp}) is zero by the same arguments as in Lemma \ref{At}.  The semiregularity map therefore descends to the absolute obstruction theory,

\[
\begin{tikzcd}[row sep=small, column sep=small]
&\BT_{\mathfrak{M}_{g,N}}[-1]  \arrow[r] &\BE^{\bullet} \arrow[r] \arrow[d,"\mathsf{sr}"] &\BE_{\mathrm{abs}}^{\bullet}  \arrow[dl, dashed, "\mathsf{sr}_{\mathrm{abs}}"] \\
&  &\CO_{M^{ \epsilon}_{\mathbf{v},\check{\beta}}(S \times C_{g,N})}[-1] &
\end{tikzcd}
\]
Hence the results of \cite{KiL} apply. 
\\

\label{chws}   The obstruction theory of the master space $MQ^{\epsilon_0}_{g,N}(M(\mathbf{v}), \beta)$  also has a surjective cosection, see \cite[Section 7.2]{N} for the definition of the master space. Like in the case of $Q^{\epsilon}_{g,N}(M(\mathbf{v}) \beta)$, it is constructed by identifying the master space  with a moduli spaces of sheaves. 

\subsection{Invariants} In what follows, we use Kiem--Li's construction of reduced virtual fundamental classes via the cosection localisation \cite{KiL}. By using the identification (\ref{idents}), we define
\[[Q^{\epsilon}_{g,N}(M(\mathbf{v}), \beta)]^{\mathrm{red}} \in H_{*}(Q^{\epsilon}_{g,N}(M(\mathbf{v}), \beta))\]
to be the associated reduced virtual fundamental class. From now on, by a virtual fundamental class we always will mean a \textit{reduced} virtual fundamental class, except for $\beta=0$, since the standard virtual fundamental class does not vanish in this case. 
These classes can be seen as virtual fundamental classes associated to the reduced obstruction-theory complex $\BE^{\bullet}_{\mathrm{red}}$, defined as the cone of the cosection, 
\begin{equation} \label{redob}
\BE_{\mathrm{red }}^{\bullet}=\mathrm{cone}(\mathsf{sr_{\mathrm{abs}}})[-1] \rightarrow \BE_{\mathrm{abs}}^{\bullet} \xrightarrow{\mathsf{sr}_{\mathrm{abs}}} \CO[-1].
\end{equation}
However, we are unable to prove that $\BE_{\mathrm{red }}^{\bullet}$ defines an obstruction theory in full generality, cf. Proposition \ref{reduction}. Luckily, Kiem--Li's class is good enough for all purposes, e.g.\ see \cite{CKL} for the virtual torus localisation of Kiem--Li's reduced classes.

\begin{definition} \label{descendent} We define  \textit{descendent quasimap} invariants, 
\[\langle \gamma_{1}\psi^{k_1}, \dots, \gamma_{N}\psi^{k_N} \rangle^{M(\mathbf{v}), \epsilon}_{g,N,\beta}:= \int_{[Q^{\epsilon}_{g,N}(M(\mathbf{v}),\beta)]^{\mathrm{red}}}\prod^{i=N}_{i=1}\ev^{*}_{i}(\gamma_{i})\psi_{i}^{k_{i}},\]
where $\gamma_{1}, \dots, \gamma_{N} \in H^{*}(M(\mathbf{v}))$, $\psi_{1}, \dots \psi_{N}$ are $\psi$-classes associated to markings, and $k_{1}, \dots k_{N}$ are non-negative integers. We similarly define the perverse invariants $\langle \gamma_{1}\psi^{k_1}, \dots, \gamma_{N}\psi^{k_N} \rangle^{\sharp,S^{[n]}, \epsilon}_{g,\beta}$, using perverse quasimaps from \cite[Section 6.3]{N}. With respect to the identification (\ref{idents}), primary $\epsilon$-invariants (no $\psi$-classes) correspond to relative DT invariants.  
\end{definition}

 Consider now the diagram 
\[
\begin{tikzcd}[row sep=small, column sep = small]
&  S \times \CC^{\epsilon}_{g,N} \arrow{dl}[swap]{p} \arrow{dr}{\pi} \\
S\times C_{g,N} && Q^{\epsilon}_{g,N}(M(\mathbf{v}), \beta)
\end{tikzcd}
\]
where $C_{g,N}$ is the universal curve over $\Mbar_{g,N}$ and $p$ is the stabilisation of curves. For the unstable values of 
$g$ and $N$, we set the product $S\times \Mbar_{g,N+1}$ to be $S$.

\begin{definition}
For a class $\tilde {\gamma} \in H^l(S\times C_{g,N})$ define the following operation on cohomology, 
\[\ch_{k+2}(\tilde{\gamma})\colon  H_{*}(Q^{\epsilon}_{g,N}(M(\mathbf{v}), \beta)) \rightarrow H_{*-2k+2-l}(Q^{\epsilon}_{g,N}(M(\mathbf{v}), \beta)),\]
\[\ch_{k+2}(\tilde{\gamma})(\xi)=\pi_{*}\left(\ch_{k+2}(\BF)\cdot p^*(\tilde{\gamma})\cap \pi^*(\xi) \right).\]
The \textit{descendent DT} invariants are then defined by
\begin{multline*} 
	\langle  \tilde{\tau}_{k_{1}}(\tilde{\gamma}_{1}), \dots, \tilde{\tau}_{k_{r}}(\tilde{\gamma}_{r}) \rangle^{\epsilon}_{g,N,\beta} \\
	=(-1)^{k_{1}}\ch_{k_{1}+2}(\tilde{\gamma}_{1}) \ \circ \ \dots \ \circ \ (-1)^{k_{r}} \ch_{k_{r}+2}(\tilde{\gamma}_{r}) \left([Q^{\epsilon}_{g,N}(M(\mathbf{v}), \beta)]^{\mathrm{red}}\right).
\end{multline*}
By combing Definition \ref{descendent} with the definition above, we  obtain a mix of descendent quasimap invariants and descendent DT  invariants, 
\[\langle  \gamma_{1}\psi^{k_1}, \dots, \gamma_{N}\psi^{k_N} \mid  \tilde{\tau}_{k_{1}}(\tilde{\gamma}_{1}), \dots, \tilde{\tau}_{k_{r}}(\tilde{\gamma}_{r}) \rangle^{M(\mathbf{v}),\epsilon}_{g,N,\beta}.\]
The same applies to the perverse invariants. For $\epsilon=0^+$ and for a fixed marked curve\footnote{We take a fiber of $Q^{0^+}_{g,N}(M(\mathbf{v}), \beta)$ over a non-stacky closed point $[(C,\mathbf{p})]\in \FM_{g,N}(\BC)$.} $(C,\mathbf{p})$, we denote the invariants above by
\[\langle  \gamma_{1}\psi^{k_1}, \dots, \gamma_{N}\psi^{k_N} \mid  \tilde{\tau}_{k_{1}}(\tilde{\gamma}_{1}), \dots, \tilde{\tau}_{k_{r}}(\tilde{\gamma}_{r}) \rangle^{S\times C}_{\check{\beta}}, \]
to emphasize that we are considering the DT theory of $S\times C$. 
\end{definition}

\subsection{Wall-crossing}

We fix a parametrized projective line $\p^1$ with a $\BC^*$-action, 
\[t[x:y]=[tx:y], \quad t\in \BC^{*},\]
such that $0:=[0:1]$ and $\infty:=[1:0]$. By convention, we set 
\[z:=e_{\BC^*}(\BC_{\mathrm{std}}),\] 
where $\BC_{\mathrm{std}}$ is the weight 1 representation of $\BC^*$. 
We now define a \textit{Vertex space}, 
\[V(M(\mathbf{v}),\beta),\]
to be a moduli space of sheaves $F$ on $S\times \p^1$ subject to the following conditions:
\begin{itemize}
\item $F$ is torsion-free, 
\item $\ch(F)=(\ch(\mathbf{v}),\check{\beta})$, 
\item fibers $F_p$ are stable for a general $p \in \p^1$, 
\item the fiber $F_\infty$  is stable, 
\item $\det(p_{C*}(p_{S}^{*}\mathbf{u}\otimes F))\cong \CO_{C}$.
\end{itemize}
By construction, there is a natural evaluation map, 
\begin{align*}\ev\colon V(M(\mathbf{v}),\beta) &\rightarrow M(\mathbf{v}),\\
F&\mapsto F_\infty.
\end{align*}
Moreover, by acting on $\p^1$, we obtain a $\BC^*$-action on $V(M(\mathbf{v}),\beta)$.  By the same arguments as in Section \ref{surjcos}, the  obstruction theory of  $V(M(\mathbf{v}),\beta)$ has a surjective cosection.  Moreover, the cosection is $\BC^*$-equivariant.
The Vertex space is not proper, but its $\BC^*$-fixed locus 
\[V(M(\mathbf{v}),\beta)^{\BC^*}\] is proper. Indeed, this follows from properness of the space of all prestable quasimaps from $\p^1$ and the fact that $V(M(\mathbf{v}),\beta)^{\BC^*}$ is just a connected  component of its $\BC^*$-fixed locus. We can therefore use the virtual torus localisation of reduced classes \cite{CKL} to define its virtual fundamental class,
\[[V(M(\mathbf{v}),\beta)]^{\mathrm{red}}:=\frac{[V(M(\mathbf{v}),\beta)^{\BC^*}]^{\mathrm{red}}}{e_{\BC^*}(\CN^{\mathrm{vir}})} \in H_*(V(M(\mathbf{v}),\beta)^{\BC^*})[z^{\pm}],\]
where $\CN^{\mathrm{vir}}$ is the virtual normal complex of $V(M(\mathbf{v}),\beta)^{\BC^*}$ inside $V(M(\mathbf{v}),\beta)$, and $z$ is the equivariant parameter. We are now ready to define Givental's $I$-function, introduced in the context of GIT quasimaps in \cite{CFKM}.

\begin{definition}We define 
\begin{align*}
	I^{\mathbf{v}}_\beta(z)&:= \ev_* [V(M(\mathbf{v}),\beta)]^{\mathrm{red}} \in H^*(M(\mathbf{v}))[z^\pm], \\
	\mu^{\mathbf{v}}_\beta(z)&:=[zI^{\mathbf{v}}_\beta(z)]_{z^{\geq 0}}  \in H^*(M(\mathbf{v}))[z].
\end{align*}

\end{definition}

\begin{lemma} \label{mu} The class $\mu^{\mathbf{v}}_\beta(z)$ admits the following expression,
\begin{equation*}
	\mu^{\mathbf{v}}_{\beta}(z)=\mu^{\mathbf{v}}_\beta \cdot\mathbbm{1}\in H^{2\dim}(M(\mathbf{v})),
\end{equation*}
where $\mu^{\mathbf{v}}_\beta \in \BQ$.
\end{lemma}

\textit{Proof.} This follows from the definition and the fact that the reduced virtual dimension of $V(M(\mathbf{v}),\beta)$ is equal to $\dim(M(\mathbf{v}))+1$. 
\qed 
\\

The space $\BR_{>0} \cup \{0^+,\infty\}$ of $\epsilon$-stabilities is divided into chambers, in which the moduli space $Q^{\epsilon}_{g,N}(M(\mathbf{v}),\beta)$ stays the same, and as $\epsilon$ crosses the a wall between chambers, the moduli space changes discontinuously. Given a quasimap class $\beta \in \Eff(M(\mathbf{v}), \rCoh(S,\mathbf{v}))$, then for a class $\beta'\in \Eff(M(\mathbf{v}), \rCoh(S,\mathbf{v}))$ that appears as a summand of $\beta'$, we define 
\[\deg(\beta'):=\beta'(\CL_\beta),\]
and we define the $\epsilon$-stability of quasimaps of degree $\beta'$ with respect to the line bundle $\CL_\beta$, constructed in \cite[Section 3.4]{N}. Let $\epsilon_{0}=1/d_{0} \in \BR_{>0}$ be a wall for $\beta \in \Eff(M(\mathbf{v}), \rCoh(S,\mathbf{v}))$, and $\epsilon_-$, $\epsilon_+$ be some values that are close to $\epsilon_{0}$ from the left and the right of the wall, respectively. 
\begin{theorem}\label{reducedws}
Assuming $2g-2+N+\epsilon_0\deg(\beta)>0$, we have 
\begin{multline*}\langle  \gamma_{1}\psi^{k_1}, \dots, \gamma_{N}\psi^{k_N}  \rangle^{M(\mathbf{v}),\epsilon_{-}}_{g,N,\beta}-\langle  \gamma_{1}\psi^{k_1}, \dots, \gamma_{N}\psi^{k_N} \rangle^{M(\mathbf{v}),\epsilon_{+}}_{g,N,\beta}\\ =\mu^{\mathbf{v}}_\beta \cdot \langle  \gamma_{1}\psi^{k_1}, \dots, \gamma_{N}\psi^{k_N} , \mathbbm{1} \rangle^{M(\mathbf{v}),\infty}_{g,N+1,0}
\end{multline*}
if $\deg(\beta)=d_{0}$, and
\[\langle  \gamma_{1}\psi^{k_1}, \dots, \gamma_{N}\psi^{k_N} \rangle^{M(\mathbf{v}),\epsilon_{-}}_{g,N,\beta}=\langle \gamma_{1}\psi^{k_1}, \dots, \gamma_{N}\psi^{k_N} \rangle^{M(\mathbf{v}),\epsilon_{+}}_{g,N,\beta},\]
otherwise. 
\end{theorem}
\textit{Sketch of Proof.} As in the case of \cite[Theorem 7.5]{N}, the results from \cite[Section 6]{YZ} apply almost without change. The difference is that we use reduced classes now. Up to a finite gerbe, the fixed components of the master space which contribute to the wall-crossing formula are of the following form,  
\[\widetilde{Q}^{\epsilon_+}_{g,N+k}(M(\mathbf{v}),\beta') \times_{M(\mathbf{v})^{k}} \prod_{i=1}^k V(M(\mathbf{v}),\beta_i)^{\BC^*},\]
where $\beta=\beta'+\beta_{1}+\dots+\beta_{k}$, and $\deg(\beta_i)=d_{0}$. Recall that  $\widetilde{Q}^{\epsilon_+}_{g,N+k}(M(\mathbf{v})\beta')$ is the base change of $Q^{\epsilon_+}_{g,N}(M(\mathbf{v}),\beta)$  from $\FM_{g,N}$ to $\widetilde{\FM}_{g,N,d}$, where the latter is the moduli space of curves with entangled tails. The reduced class of a product splits as a product of reduced and non-reduced classes on its factors (cf. \cite[Section 3.9]{MPT}). Hence by Corollary \ref{cosection} and \cite{KiL}, it vanishes, unless $\beta'=0$ and $k=1$. In this case, \[\widetilde{Q}^{\epsilon^+}_{g,N+1}(M(\mathbf{v}),0)={Q}^{\infty}_{g,N+1}(M(\mathbf{v}),0)=\Mbar_{g,N+1}(M(\mathbf{v}),0).\]
Using the analysis presented in \cite[Section 7]{YZ} and Lemma \ref{mu}, we get that the contribution of this component to the wall-crossing is
\[ \langle \gamma_{1}\psi^{k_1}, \dots, \gamma_{N}\psi^{k_N}, \mu^{\mathbf{v}}_\beta(-\psi_{N+1}) \rangle^{M(\mathbf{v}),\infty}_{g,N+1,0}= \mu^{\mathbf{v}}_\beta \cdot  \langle \gamma_{1}\psi^{k_1}, \dots, \gamma_{N}\psi^{k_N} , \mathbbm{1} \rangle^{M(\mathbf{v}),\infty}_{g,N+1,0},\]
this concludes the argument. 
\qed
\\

Applying Theorem \ref{reducedws} inductively to all walls on the way from $\epsilon=0^+$ to $\epsilon=\infty$, we obtain the following result.
\begin{corollary}	\label{cor1}Assuming $(g,N)\neq(0,1)$, we have 
\begin{multline*}\langle  \gamma_{1}\psi^{k_1}, \dots, \gamma_{N}\psi^{k_N}  \rangle^{M(\mathbf{v}),0^+}_{g,N,\beta}-\langle  \gamma_{1}\psi^{k_1}, \dots, \gamma_{N}\psi^{k_N} \rangle^{M(\mathbf{v}),\infty}_{g,N,\beta}\\
	=\mu^{\mathbf{v}}_\beta \cdot \langle  \gamma_{1}\psi^{k_1}, \dots, \gamma_{N}\psi^{k_N} , \mathbbm{1} \rangle^{M(\mathbf{v}),\infty}_{g,N+1,0}.
\end{multline*}
\end{corollary} 
Another immediate corollary of the wall-crossing formula is the following expression for the wall-crossing invariants. 
\begin{corollary} \label{cormu} We have
\begin{equation*} 
	\mu^{\mathbf{v}}_{\beta}=\langle[\mathrm{pt}],\mathbbm{1}\rangle^{M(\mathbf{v}), 0^+}_{0,2,\beta}.
\end{equation*}
\end{corollary}
\textit{Proof.} The result follows from Corollary \ref{cor1}  applied to the invariants $\langle [\mathrm{pt}], \mathbbm{1}, \rangle^{M(\mathbf{v}),\epsilon}_{0,3,\beta}$ for $\epsilon \in \{0^+, \infty\}$, and the string equation on the GW side. 
\qed 
\\

There are invariants that are not covered by the results above and of great interest for us - those of a fixed elliptic curve.  Let $E$ be a fixed elliptic curve, and  $Q^\epsilon_{E}(M(\mathbf{v}),\beta)^{\bullet}$ be the fiber of 
\[Q^\epsilon_{1,0}(M(\mathbf{v}),\beta) \rightarrow \Mbar_{1,0}\]
over the stacky point $[E]/E \in  \Mbar_{1,0}(\BC)$. In other words, 
$Q^\epsilon_{E}(M(\mathbf{v}),\beta)^{\bullet}$ is the moduli space of $\epsilon$-stable quasimaps, whose smoothing of the domain is $E$, and maps are considered up translations of $E$. For $\beta\neq 0$, we define 
\[ \langle \emptyset \rangle^{M(\mathbf{v}),\epsilon}_{E,\beta}:=\int_{[Q^{\epsilon}_{E}(M(\mathbf{v}),\beta)^{\bullet}]^{\mathrm{red}}}1.\]

\begin{theorem} \label{elliptic} Assuming $\beta\neq0$, we have
\[ \langle \emptyset \rangle^{M(\mathbf{v}),\epsilon_-}_{E,\beta}= \langle \emptyset \rangle^{M(\mathbf{v}),\epsilon_+}_{E,\beta}+ \mu^{\mathbf{v}}_{\beta}\cdot \chi(M(\mathbf{v}))  ,\]
if $\deg(\beta)=d_{0}$, and 
\[\langle \emptyset \rangle^{M(\mathbf{v}),\epsilon_-}_{E,\beta}= \langle \emptyset \rangle^{M(\mathbf{v}),\epsilon_+}_{E,\beta},\]
otherwise. 

\end{theorem}
\textit{Sketch of Proof.} As in Theorem \ref{reducedws}, the only case when the contribution from the wall-crossing components is non-zero is the one of $\beta'=0$ and $k=1$. In this case, \[\widetilde{Q}^{\epsilon^+}_{(E,0_E)}(M(\mathbf{v}),0)\cong M(\mathbf{v}),\] 
and the obstruction bundle is the tangent bundle $T_{M(\mathbf{v})}$. Hence the virtual fundamental class is $\chi(M(\mathbf{v}))[\mathrm{pt}]$. The corresponding wall-crossing term is therefore equal to
\[\mu^{\mathbf{v}}_\beta \cdot \chi(M(\mathbf{v})),\]
this concludes the argument. 
\qed
\\

As before, by applying Theorem \ref{elliptic} inductively, we obtain the following result.
\begin{corollary}  \label{cor2} Assuming $\beta\neq 0$, we have 
\[ \langle \emptyset \rangle^{M(\mathbf{v}),0^+}_{E,\beta}= \langle \emptyset \rangle^{M(\mathbf{v}),\infty}_{E,\beta}+ \mu^{\mathbf{v}}_\beta \cdot  \chi(M(\mathbf{v})). \]
\end{corollary}

\section{Applications}

\subsection{Genus-0 invariants of $S^{[n]}$} 
We start with genus-0 3-point quasimap invariants of $S^{[n]}$.   
The moduli space of $\p^1$ with  3 marked points is  a point. Hence by fixing markings, we can identify moduli spaces of $0^+$-stable quasimaps with 1-dimensional subschemes/stable pairs on $S\times \p^1$ relative to the divisor $S_{0,1,\infty} \subset S\times \p^1$. In the notation introduced in (\ref{not}) and (\ref{idents}), we therefore obtain
\begin{align} \label{ident}
\begin{split}
	Q^{0^+}_{0,3}(S^{[n]}, \beta)&\cong \mathrm{Hilb}^{0^+}_{n,\check{\beta}}(S\times C_{0,3}) =\mathrm{Hilb}_{n,\check{\beta}}(S\times \p^1/S_{0,1, \infty}),\\
	Q^{0^+}_{0,3}(S^{[n]}, \beta)^\sharp&\cong \mathrm{P}^{0^+}_{n,\check{\beta}}(S\times C_{0,3}) =\mathrm{P}_{n,\check{\beta}}(S\times \p^1/S_{0,1, \infty}),
\end{split}
\end{align}
such that relative insertions correspond to primary quasimap insertions.
Moreover, by Corollary \ref{cor1} and the string equation, the wall-crossing is trivial for primary invariants, if $(g,N)=(0,3)$. We  therefore obtain that 
\[\langle\gamma_{1},\gamma_{2},\gamma_{3}\rangle_{0,3,\beta}^{ S^{[n]},0^+}=\langle\gamma_{1},\gamma_{2},\gamma_{3}\rangle_{0,3,\beta}^{S^{[n]},\infty}=\langle\gamma_{1},\gamma_{2},\gamma_{3}\rangle_{0,3,\beta}^{\sharp, S^{[n]},0^+}.\]
In light of the identification (\ref{ident}), we obtain  the following result.
\begin{corollary} \label{cor3} We have
\[\langle\gamma_{1},\gamma_{2},\gamma_{3}\rangle_{n,\check{\beta}}^{ S\times \p^1}=\langle\gamma_{1},\gamma_{2},\gamma_{3}\rangle_{0,3,\beta}^{S^{[n]},\infty}=\langle\gamma_{1},\gamma_{2},\gamma_{3}\rangle_{n,\check{\beta}}^{\sharp, S\times \p^1}.\]
\end{corollary}

On one hand,  the result above together with the PT/GW correspondence for $K3$ geometries of \cite[Theorem 1.2]{Ob21} confirm the conjecture proposed in \cite[Conjecture 1]{OP1}. The conjecture claimed that the relative GW theory of $S\times \p^1$ and the genus-0 3-point GW theory of $S^{[n]}$ are equivalent after the change of variables $y=-e^{iu}$. In fact, Corollary \ref{cor3} is a more natural form of \cite[Conjecture 1]{OP1}, as it asserts equality of invariants without any change of variables. On the other hand, Corollary \ref{cor3} provides a DT/PT correspondence with relative insertions for $S\times \p^1$.  More generally, the results above can be restated for a relative geometry
\[S\times C_{g,N} \rightarrow \Mbar_{g,N},\] 
such that $N>2$. In this case, by the string equation, the wall-crossing is also trivial for primary insertions.

\subsection{Genus-1 invariants of $S^{[n]}$and Igusa cusp form conjecture} \label{Sectiongenus1} In this section, we will consider perverse quasimaps. In this case, $0^+$-stable quasimaps correspond to stable pairs. Let us firstly establish a relation between degrees $\beta$ of quasimaps and Chern characters $\check{\beta}$ of the associated stable pairs. 

Firstly, the homology  $H_2(S^{[n]},\BZ)$ admits a Nakajima basis
\begin{align*}
H_2(S,\BZ) \oplus \BZ &\xrightarrow{\sim} H_2(S^{[n]},\BZ), \\
(\gamma, k) &\mapsto C_{\gamma}+kA,
\end{align*}
where the classes above are defined in terms of Nakajima operators as follows, 
\[C_\gamma=\mathfrak{q}_{-1}(\gamma)\mathfrak{q}_{-1}([\mathrm{pt}])^{n-1}1_S, \quad  A=\mathfrak{q}_{-2}([\mathrm{pt}])\mathfrak{q}_{-1}([\mathrm{pt}])^{n-2}1_S.\]
We refer to \cite[Section 1]{Ob18} for the notation and the definition of Nakajima operators in the similar context. 
In more geometric terms, if a class $\gamma$ is represented by a curve $\Gamma \subset S$, then the class $C_\gamma$ is represented by the curve $\Gamma_n \subset S^{[n]}$ which is given by letting one point move along $\Gamma$ and keeping $n-1$ other distinct points fixed. The class $A$ is given by the locus of  length-2 non-reduced structures on a fixed point $p\in S$, keeping other $n-2$ reduced points fixed.

Given now a quasimap $f \colon C \rightarrow \rpCoh(S,\mathbf{v})$.  By \cite[Theorem 6.8]{N}, the objects associated to such quasimaps are stable pairs on $S\times C$, hence the class $-\check{\beta}$ is of the following form 

\begin{equation} \label{cohidentification}
-\check{\beta}=(0,\gamma, k)\in H^{0}(S)\oplus H^{1,1}(S)\oplus H^{4}(S),
\end{equation}
where the negative sign amounts to considering Chern characters of subschemes rather than their ideals. For simplicity we will denote a class $-\check{\beta}$ just by $(\gamma, k)$. The decomposition above and the one given by Nakajima basis are related. If we consider a map $f \colon C \rightarrow S^{[n]}$ of degree $\beta$ (in the sense of quasimaps), such that  $-\check{\beta}=(\gamma, k)$, then by \cite[Lemma 2]{Ob18} we have
\begin{equation} \label{nakbeta}
f_*[C]=C_\gamma +kA. 
\end{equation}
Hence from now on, we will write degrees of quasimaps $\beta$ in terms of $-\check{\beta}$, which by (\ref{nakbeta}) also corresponds to writing degrees in terms Nakajima basis in the case $\epsilon=\infty$. 

Consider now a generic projective $K3$ surface $S$ with of Picard rank 1, such that
\[\mathrm{NS}(S)=\BZ \langle \beta_{h} \rangle ,  \quad \beta_{h}^2=2h-2.\] 
By  the previous discussion and \cite[Theorem 6.8]{N}, we have the following identification of moduli spaces
\[Q^{0^+}_{E}(S^{[n]},(\beta_{h},k))^{\bullet,\sharp}\cong[\mathrm{P}_{n,(\beta_{h},k)}(S\times E)/E].\]
As before, the superscript on the moduli space on the left indicates that we consider maps up to translations of $E$. For the same reason, we take the quotient by $E$ on the left. On the other hand,
\[Q^{\infty}_{E}(S^{[n]},(\beta_{h},k))^{\bullet,\sharp}=\Mbar_{E}(S^{[n]},(\beta_{h},k))^{\bullet}.\] Consider now the following two generating series, 
\begin{align*}
\mathsf{PT}(p,q,\tilde{q})&:=\sum_{n\geq0}\sum_{h\geq0}\sum_{k\in \BZ}(-p)^{k}q^{h-1}\tilde{q}^{n-1} \langle \emptyset \rangle^{S^{[n]},0^+}_{E,(\beta_{h},k)},\\
\mathsf{GW}(p,q,\tilde{q})&:=\sum_{n>0}\sum_{h\geq0}\sum_{k\geq0} (-p)^{k}q^{h-1}\tilde{q}^{n-1}\langle \emptyset \rangle^{S^{[n]},\infty}_{E,(\beta_{h},k)}.
\end{align*}
The series are well-defined, because $(S,\beta)$ and $(S',\beta')$ are deformation equivalent, if and only if \[\beta^{2}=\beta'^2 \quad \text{and} \quad \mathrm{div}(\beta)=\mathrm{div}(\beta'),\] 
where $\mathrm{div}(\beta)$ is the divisibility of the class. In our case, the classes $\beta_{h}$ are primitive by definition. In \cite{OS2} and \cite{OPi}, it was proved that
\[\mathsf{PT}(p,q,\tilde{q})=\frac{1}{-\chi_{10}(p,q,\tilde{q})},\]
where $\chi_{10}(p,q,\tilde{q})$ is the \textit{Igusa cusp form}, we refer to \cite[Section 0.2]{OPi} for its definition.  We can view both series as generating series of quasimaps for $\epsilon\in \{0^+,\infty\}$. Using Corollary \ref{cor2}, we obtain 
\begin{equation*}\mathsf{PT}(p,q,\tilde{q})=\mathsf{GW}(p,q,\tilde{q}) 
+\sum_{n\geq0}\sum_{h\geq0}\sum_{k\in \BZ} (-p)^{k}q^{h-1}\tilde{q}^{n-1} \mu^{n,\sharp}_{(\beta_{h},k)} \cdot \chi(S^{[n]}).
\end{equation*}
The wall-crossing invariants  $\mu^{n,\sharp}_{(\beta_{h},k)}$
are equal to virtual Euler characteristics of Quot schemes, as it is explained in \cite{Ob}.  In the same article,  wall-crossing invariants are also explicitly computed for $S^{[n]}$, see \cite[Theorem 1.2]{Ob}, 
\[\sum_{n\geq0}\sum_{h\geq0}\sum_{k\in \BZ} (-p)^{k}q^{h-1}\tilde{q}^{n-1} \mu^{n,\sharp}_{(\beta_{h},k)}  \chi(S^{[n]})=\frac{1}{\Theta^2\Delta}\frac{1}{\tilde{q}}\prod_{n\geq1}\frac{1}{(1-(\tilde{q} G)^n)^{24}},\]
where 
\begin{align*}
\Theta(p,q)&=(p^{1/2}-p^{1/2})\prod_{m\geq1} \frac{(1-pq^m)(1-p^{-1}q^m)}{(1-q^m)^2} ,\\
G(p,q)&=-\Theta(p,q)^2 \left(p\frac{d}{dp}\right)^2\log(\Theta(p,q)),
\end{align*}
and $\Delta(q)=q\prod_{n\geq 1}(1-q^n)^{24}.$ We therefore obtain the following corollary, which confirms the first equality in \cite[Conjecture A]{OPa}. 
\begin{corollary} We have
\[\mathsf{PT}(p,q,\tilde{q})=\mathsf{GW}(p,q,\tilde{q})+ 
\frac{1}{\Theta^2\Delta}\frac{1}{\tilde{q}}\prod_{n\geq1}\frac{1}{(1-(\tilde{q} G)^n)^{24}}.\]

\end{corollary}
As in the case of \cite[Conjecture 1]{OP1}, the statement of \cite[Conjecture A]{OPa} involved the relative GW theory of $S\times E$ instead of the relative DT theory of $S\times E$. In light of our wall-crossing, the latter should be considered more natural in this context. Nevertheless, the two are equivalent by \cite{Ob21}. 

\subsection{Higher-rank DT invariants}\label{higher} Under our assumptions, a moduli space $M(\mathbf{v})$ is  deformation equivalent to a punctorial Hilbert scheme $S^{[n]}$, where $2n=\dim(M(\mathbf{v}))$.  Since we are working with reduced classes, in order to use the deformation invariance of GW theory,  we have to deform $M(\mathbf{v})$ together with a curve class $\beta$. In this case, we might need to deform  $(M(\mathbf{v}),\beta)$ to a Hilbert scheme of points  $(S'^{[n]},\beta')$ on another $K3$ surface $S'$. We conclude that the reduced GW theory of $M(\mathbf{v})$ in the class $\beta$ is equivalent to the one of $S'^{[n]}$ in the class $\beta'$. Applying the quasimap wall-crossing both to $M(\mathbf{v})$ and to $S'^{[n]}$, we can therefore express higher-rank DT invariants of a threefold $S\times C$ in terms of rank-one DT invariants and wall-crossing invariants.

\subsubsection{$K3\times \p^1$}  \label{p1}
Let us firstly consider invariants on $S\times \p^1$ relative to $S_{0,1,\infty}$. As previously, by fixing the markings, we obtain
\[M^{ 0^+}_{\mathbf{v},\check{\beta}}(S \times C_{0,3})=M_{\mathbf{v},\check \beta}(S\times \p^1/S_{0,1, \infty}).\]
The $\BC$-valued points of $M_{\mathbf{v},\check \beta}(S\times \p^1/S_{0,1, \infty})$ are triples  $(P,\{0,1,\infty\}, F)$ satisfying the same conditions as triples in (\ref{triples}), such that:
\begin{itemize}
\item $(P, \{0,1,\infty\})$ is a  marked curve without rational tails\footnote{Rational components with one special point, i.e.\ with one separating marking or one node.}, whose smoothing\footnote{Equivalently, $P$ is  an isotrivial  degeneration  of $\p^1$ at $\{0,1,\infty\}$, i.e. a bubbling of $\p^1$ at $\{0,1,\infty\}$.} is $(\p^1, \{0,1,\infty\})$.
\end{itemize}
Moreover, as in the case of $S^{[n]}$, there is no wall-crossing  by Theorem \ref{cor1} and the string equation, therefore
\[\langle\gamma_{1},\gamma_{2},\gamma_{3}\rangle_{0,3,\beta}^{M(\mathbf{v}),0^+}=\langle\gamma_{1},\gamma_{2},\gamma_{3}\rangle_{0,3,\beta}^{M(\mathbf{v}),\infty}.\]
Choose a deformation of $(M(\mathbf{v}),\beta)$ to $(S'^{[n]},\beta')$ which keeps  the curve class $\beta$ algebraic. The deformation gives an identification of cohomologies  
\[ H^{*}(M(\mathbf{v}))\cong H^*(S'^{[n]}),\] 
which we use to identify curve classes $\beta$ and $\beta'$, and insertions $\gamma_1,\dots, \gamma_N$ on both sides. 
With respect to this identification, we have
\begin{equation} \label{equality}
\langle\gamma_{1},\gamma_{2},\gamma_{3}\rangle_{0,3,\beta}^{M(\mathbf{v}),0^+}=\langle\gamma_{1},\gamma_{2},\gamma_{3}\rangle_{0,3,\beta}^{M(\mathbf{v}),\infty}
=\langle\gamma_{1},\gamma_{2},\gamma_{3}\rangle_{0,3,\beta}^{S'^{[n]},\infty}=\langle\gamma_{1},\gamma_{2},\gamma_{3}\rangle_{0,3,\beta}^{S'^{[n]},0^+}.
\end{equation}
Passing from quasimaps to sheaves and using $(\ref{equality})$, we obtain the following result. 
\begin{corollary} Given a deformation of $(M(\mathbf{v}),\beta)$ to $(S'^{[n]},\beta')$. Identifying cohomologies of  $M(\mathbf{v})$ and  $S'^{[n]}$, we have 

\[\langle\gamma_{1},\gamma_{2},\gamma_{3}\rangle_{ \mathbf{v},\check{\beta}}^{S\times \p^1}=\langle\gamma_{1},\gamma_{2},\gamma_{3}\rangle_{n,\check{\beta}}^{S'\times \p^1}.\]
\end{corollary}

\subsubsection{$K3\times E$} \label{ell} Consider now $S\times E$, where $E$ is an elliptic curve. If $M(\mathbf{v})=S^{[n]}$, then invariants  $\langle[\mathrm{pt}],\mathbbm{1}\rangle^{ S^{[n]},0^+}_{0,2,\beta}$ from Corollary \ref{cormu} are called  \textit{rubber} DT  invariants on $S\times \p^1$. These are invariants associated to the moduli space of subschemes on $S\times \p^1$ relative to the divisor $S_{0,\infty}$ 
up to the $\BC^{*}$-action coming from $\p^1$-factor which fixes $0$ and $\infty$, 
\[[\mathrm{Hilb}_{n,\check{\beta}}(S\times \p^1/S_{0, \infty})/\BC^*].\]
These invariants can be \textit{rigidified} to standard relative DT invariants with absolute insertions. In fact, this holds more generally for any $\mathbf{v}$ and also for genus-1 invariants $\langle \emptyset \rangle^{M(\mathbf{v}),0^+}_{E,\beta}$.

\begin{lemma} \label{rig} We have
\begin{align*}
	\langle  \tilde{\tau}_0( D\boxtimes \omega) \rangle^{S\times E}_{\mathbf{v},\check{\beta}}&=( \check{\beta}_1\cdot D) \langle \emptyset \rangle^{M(\mathbf{v}),0^+}_{E,\beta},  \\
	\langle [\mathrm{pt}], \mathbbm{1} \mid \tilde{\tau}_0\langle D\boxtimes \omega\rangle \rangle^{S\times \p^1}_{\mathbf{v},\check{\beta}}&=(\check{\beta}_1\cdot D)\langle[\mathrm{pt}],\mathbbm{1}\rangle^{M(\mathbf{v}),0^+}_{0,2,\beta} ,
\end{align*}
where $D \in H^2(S)$, $\omega \in H^2(C)$ is the point class, and $\check{\beta}_1 \in H^2(S)$ is the component of $\check{\beta}$ of cohomological degree 2.
\end{lemma}
\textit{Proof.}  The proof is exactly the same as in \cite[Lemma 3.3]{MO2}. For  $S\times E$, see also \cite[Section 3.4]{OberG}  \qed
\\

Applying the same procedure as for $S\times \p^1$, using Corollary \ref{cormu} and  Lemma \ref{rig}, we obtain the following result.

\begin{corollary}  We have
\begin{multline*}
	\langle  \tilde{\tau}_0( D\boxtimes \omega) \rangle^{S\times E}_{\mathbf{v},\check{\beta}}=	\langle  \tilde{\tau}_0( D\boxtimes \omega) \rangle^{S'\times E}_{n,\check{\beta}} \\
	+\chi(S^{[n]})	\left(\langle [\mathrm{pt}], \mathbbm{1} \mid \tilde{\tau}_0( D\boxtimes \omega)\rangle^{S\times \p^1}_{\mathbf{v},\check{\beta}}- 	\langle [\mathrm{pt}], \mathbbm{1} \mid \tilde{\tau}_0( D\boxtimes \omega)\rangle^{S'\times \p^1}_{n,\check{\beta}}\right).
\end{multline*}
\end{corollary}

By degenerating $\p^1$ to $\p^1 \cup \p^1$, sending the interior marking and the relative marking to the first component, and applying the degeneration formula, we obtain
\[\langle [\mathrm{pt}], \mathbbm{1} \mid \tilde{\tau}_0(D\boxtimes \omega)  \rangle^{S\times \p^1}_{\mathbf{v},\check{\beta}}=\langle [\mathrm{pt}]\mid \tilde{\tau}_0(D\boxtimes \omega) \rangle^{S\times \p^1}_{\mathbf{v},\check{\beta}},\]
where we used the fact that a reduced class restricts to a reduced and non-reduced classes  on irreducible components,  which implies that it vanishes, unless $\check{\beta}=0$ on one of the components.  We refer to \cite[Section 3.4]{MNOP} (see also \cite{LW}) for the standard  degeneration formula, and to \cite[Section 5.1]{Ob}  for the reduced one. 

Similarly, by degenerating $E$ to $E \cup \p^1$, sending the interior marking to the second component, and applying the degeneration formula, we obtain 
\begin{align*}
\langle  \tilde{\tau}_0(D\boxtimes \omega)  \rangle^{S\times E}_{\mathbf{v},\check{\beta}}=  \langle  \mathbbm{1} \mid\tilde{\tau}_0(D\boxtimes \omega)   \rangle^{S\times E}_{\mathbf{v},\check{\beta}}+\chi(M(\mathbf{v}))\langle [\mathrm{pt}]\mid  \tilde{\tau}_0(D\boxtimes \omega)  \rangle^{S\times \p^1}_{\mathbf{v},\check{\beta}}.
\end{align*}
The second term on the right is the wall-crossing term. Hence we obtain the following equality of DT invariants on $S\times E$. 
\begin{corollary}  We have
\[\langle  \mathbbm{1} \mid  \tilde{\tau}_0(D\boxtimes \omega)  \rangle^{S\times E}_{\mathbf{v},\check{\beta}}=\langle  \mathbbm{1} \mid \tilde{\tau}_0(D\boxtimes \omega)   \rangle^{S'\times E}_{n,\check{\beta}}.\]
\end{corollary} 
Using the Igusa cusp form conjecture, we can obtain an explicit expression for these higher-rank relative DT invariants. Moreover, by \cite[Lemma 4.13]{N}, the higher-rank invariants associated to the moduli space $M_{\mathbf{v},\check{\beta}}(S\times E)$ can be related to invariants associated to moduli spaces of sheaves with a fixed determinant, denoted by $\widetilde{M}_{\mathbf{v},\check{\beta}}(S\times E)$,
\[\int_{[\widetilde{M}_{\mathbf{v},\check{\beta}}(S\times E)/E]^{\mathrm{red}}}1=\rk(\mathbf{v})^2 \langle \emptyset \rangle^{M(\mathbf{v}),0^+}_{E,\beta},\]
while the stability of fibers can be related to the slope stability on the threefold by \cite[Corollary A.7]{N}.


\appendix
\section{Reduced obstruction theory} \label{reduced} Consider the reduced obstruction-theory complex  $\BE_{\mathrm{red}}^{\bullet}$ from (\ref{redob}). In this section, under certain assumptions, we will construct the obstruction-theory morphism, \[(\BE_{\mathrm{red}}^{\bullet})^{\vee} \rightarrow  \BL_{Q^{\epsilon}_{g,N}(M(\mathbf{v}), \beta) / \FM_{g,N}}.\]
The proof closely follows \cite{KT}.  

\begin{proposition} \label{reduction}
Given $(\ch(\mathbf{v}),\check{\beta}) \in \Lambda\oplus \Lambda(-1)$. Assume a first-order deformation $\kappa_S \in HT^2(S)\cong HH^2(S)$ from Proposition \ref{nonc} is represented by a $\BC[\epsilon]/\epsilon^2$-linear  admissible subcategory,
\[\CC \subseteq \mathrm{D_{perf}}(\CY),\]
where $\CY \rightarrow B=\Spec{\BC[\epsilon]/\epsilon^2}$ is flat. Then there exists an obstruction theory morphism,
\[(\BE_{\mathrm{red}}^{\bullet})^{\vee} \rightarrow  \BL_{Q^{\epsilon}_{g,N}(M(\mathbf{v}), \beta) / \FM_{g,N}}.\]
\end{proposition}

\textit{Proof.}
Firstly, by taking the central fiber, we get that
\[ \mathrm{D_{perf}}(S) \subseteq \mathrm{D_{perf}}(Y)\]is an admissible subcategory, where $Y$ is the central fiber of $\CY$. Therefore there is an isomorphism of moduli stacks 
\begin{equation} \label{Y}
\Coh(S)\cong \mathfrak{D}_{\Cohc(S)}(Y),
\end{equation}
where $\mathfrak{D}_{\Cohc(S)}(Y)$ is the moduli stack of objects on $Y$ which are contained in the subcategory $\Cohc(S)\subseteq \mathrm{D_{perf}}(Y)$. This also implies that  moduli stacks of quasimaps associated to pairs  $(M(\mathbf{v}),\Coh(S))$ and $(M(\mathbf{v}),\mathfrak{D}_{\Cohc(S)}(Y))$ are isomorphic,
\[Q^{\epsilon}_{g,N}(M(\mathbf{v}),\Coh(S),\beta) \cong Q^{\epsilon}_{g,N}(M,\mathfrak{D}_{\Cohc(S)}(Y),\beta ).\]
Let 
\[M_{S}:= M^{ \epsilon}_{\mathbf{v}, \check{\beta}}(S \times C_{g,N}) \cong M^{ \epsilon}_{\mathbf{v}, \check{\beta}}(Y \times C_{g,N})=: M_{Y}\] be the relative moduli spaces of objects which are associated to moduli spaces  $Q^{\epsilon}_{g,N}(M,\Coh(S),\beta)$ and $Q^{\epsilon}_{g,N}(M,\mathfrak{D}_{\Cohc(S)}(Y),\beta )$, defined as in \cite[Section 3.5]{N}.

Secondly, the inclusion 
\[\mathrm{D^b}(S) \hookrightarrow \mathrm{D_{perf}}(Y)\]induces a map between the Hochschild cohomologies, 
\begin{equation}
HH^2(Y) \rightarrow HH^2(S),
\end{equation}
given by restricting the natural transformation of functors,
\[\mathrm{id_{D_{perf}(Y)}} \rightarrow [2].\]
This map sends $\kappa_Y$ to $\kappa_S$ (see e.g.\ \cite[Lemma 4.6]{P}), where $\kappa_Y$ is the class associated to the deformation $\CY \rightarrow B$.
Moreover, for a complex $F \in \mathrm{D^{\mathrm{b}}}(S\times C)$, the class
\[\kappa(F) \in \Ext^2(F,F),\]which is given by applying the natural transformation associated to $\kappa \in HH^2(S)$ to $F$,
is the obstruction to deform $F$ in the direction of $\kappa$. By \cite[Proposition 5.2]{To2} and \cite{Cal},
it agrees with obstruction class given by composing the Kodaira--Spencer class with the Atiyah class,
\[\kappa(F)=\kappa \cdot \exp(-\mathrm{At}(F)),\]after applying the HKR isomorphism,
\[HH^2(S)\cong HT^2(S).\]We now identify a sheaf $F\in \Cohc(S\times C)$ with its image in $\mathrm{D_{perf}}(Y\times C)$, then the following triangle commutes
\[
\begin{tikzcd}[row sep=small, column sep = small]
& HH^2(S) \arrow[r] & \Ext^2(F,F) \\
& HH^2(Y) \arrow[u] \arrow[ur]
\end{tikzcd}
\]
Hence by the choice of $\kappa_S$, the deformation of sheaves in the class $(\mathbf{v},\check{\beta})$ viewed as complexes on $Y\times C$ is obstructed in the direction of $\kappa_{Y}$, because the obstruction class is non-zero by the construction of $\kappa_S$. 

Now, we will closely follow  \cite[Section 3.2]{KT}. By the above discussion the inclusion of the central fiber over $B$,
\[M_Y \hookrightarrow M_{\CY/B},\]
is an isomorphism.  The obstruction complexes of $M_{Y}$ and $M_{S}$ are isomorphic under the natural identifications of the moduli spaces,
\begin{equation} \label{ident2}
R\CH om_{\pi_S}(\BF_S,\BF_S)\cong R\CH om_{\pi_Y}(\BF_Y,\BF_Y),
\end{equation}
because both complexes can be defined just in terms of $\mathrm{D^b}(S)$, where $\BF_{S/Y}$ are universal families of $M_{S/Y}$, and $\pi_{S/Y}$ are natural projections.
Recall that we are interested in the trace-zero part  $R\CH om_{\pi_S}(\BF_S,\BF_S)_{0}$ of $R\CH om_{\pi_S}(\BF_S,\BF_S)$, cf. Section \ref{surjcos}. 

\begin{claim} The composition
\begin{equation}\label{obs}
	(\BE^{\bullet})^{\vee}=(R\CH om_{\pi_S}(\BF_S,\BF_S)_{0}[1])^{\vee}\rightarrow (R\CH om_{\pi_{Y}}(\BF_{Y},\BF_{Y})[1])^{\vee}  \rightarrow \BL_{M_{\CY/B}/B},
\end{equation}
where the first map is given by identification (\ref{ident2}) and the second is given by the Atiyah class on $\CY \times M_{\CY/B}$, is a perfect obstruction theory.
\end{claim}

\textit{Proof of Claim.}
For the proof of the claim, we plan to use the criteria from \cite[Theorem 4.5]{BeFa}. For any $B$-scheme $Z_{0}$, a $B$-map $Z_{0} \rightarrow M_{\CY/B}$ factors though the central fiber. Hence the $B$-structure map $Z_{0}\rightarrow B$ factors through the closed point of $B$.  Let $\CF_0$ be the sheaf associated to the map $Z_{0} \rightarrow M_{\CY/B}$.  The morphism 
\[(R\CH om_{\pi_{Y}}(\BF_{Y},\BF_{Y})[1])^{\vee} \rightarrow \BL_{M_{\CY/B}/B}\]is an obstruction theory. By \cite[Theorem 4.5]{BeFa}, to prove that (\ref{obs}) is an obstruction theory, it suffices to prove that the image of a non-zero obstruction class $\varpi (\CF_{0}) \in \Ext^2_{Y\times Z_{0}}(\CF_{0},\CF_{0} \otimes p_Y^*I) $ with respect to the map
\begin{equation}\label{obst}
\Ext^2_{Y\times Z_{0}}(\CF_{0},\CF_{0} \otimes p_Y^*I) \cong \Ext^2_{S\times Z_{0}}(\CF_{0},\CF_{0} \otimes p_S^*I) 
\rightarrow  \Ext^2_{S\times Z_{0}}(\CF_{0},\CF_{0} \otimes p_S^*I)_0
\end{equation} 
is non-zero for any square-zero $B$-extension $Z$ of $Z_0$ given by an ideal $I$, where $p_Y\colon Y\times_{B} Z_{0}=Y \times Z_0 \rightarrow Z_{0}$ and $p_S\colon S\times Z_0 \rightarrow Z_{0}$ are the natural projections. Given a square-zero $B$-extension $Z$ of $Z_{0}$, there are two possibilities:
\begin{enumerate}
 \item[1.] the $B$-structure map $Z\rightarrow B$ factors through the closed point;
\item[2.] the $B$-structure map $Z\rightarrow B$ does not factor through the closed point.
\end{enumerate} 
We will deal with them separately. 
\\

\noindent 1. In this case, the obstruction of lifting the map to $Z\rightarrow M_{\CY/B}$ coincides with the obstruction of lifting the map to $Z\rightarrow M_Y\cong M_S$, hence if $\varpi(\CF_{0})$ is non-zero, its image with respect (\ref{obst}) is non-zero.
\\

\noindent 2.  In this case, a lift to $Z \rightarrow  M_{\CY/B}$ is always obstructed, and the obstruction is already present at a single fiber of $p_{Y}$ in the following sense. By assumption there exists a section $B \rightarrow Z$ which is an immersion (we can find an open affine subscheme $U\subset Z$ such that $U \rightarrow B$ is flat, but then $U\cong U_{0} \times B$, because first-order deformations of affine schemes are trivial, thereby we get a section). Let $z\in Z$ be image of the closed point of $B$ of the section, then the restriction 
\[\Ext^2_{Y\times S_{0}}(\CF_{0},\CF_{0} \otimes p_{Y}^*I) \rightarrow \Ext^2_{Y\times z}( \CF_{0,z},\CF_{0,z}\otimes  p_{Y}^*I_{z})\]applied to the obstruction class $ \varpi (\CF_{0})$ is non-zero and is the obstruction to lift the sheaf $\CF_{0,z}$ on $Y$ to a sheaf on $\CY$, hence due to the following commutative diagram,

\[
\begin{tikzcd}[row sep=small, column sep = small]
& \Ext^2_{Y\times Z_{0}}(\CF_{0},\CF_{0} \otimes p_Y^*I) \arrow[d] \arrow[r] & \Ext^2_{Y\times z}( \CF_{0,z},\CF_{0,z}\otimes  p_{Y}^*I_{z}) \arrow[d] \\
&\Ext^2_{S\times Z_{0}}(\CF_{0},\CF_{0} \otimes p_S^*I)_{0} \arrow[r] & \Ext^2_{S\times z}( \CF_{0,z},\CF_{0,z}\otimes  p_S^*I_{z})_{0} 
\end{tikzcd}
\]
we conclude that the image of $\varpi (\CF_{0})$ in $\Ext^2_{S\times Z_{0}}(\CF_{0},\CF_{0} \otimes p^*I)_{0}$ is non-zero, because the image of $\varpi (\CF_{0,z})$ is non-zero in $\Ext^2_{S\times z}( \CF_{0,z},\CF_{0,z}\otimes  p^*I_{z})_{0}$. This establishes claim. \\

The absolute perfect obstruction theory $(\BH^{\bullet})^{\vee}$ is then defined by taking the cone of $ (\BE^{\bullet})^{\vee} \rightarrow \Omega_{B}[1]$. Hence we have the following diagram 
\[
\begin{tikzcd}[row sep=small, column sep = small]
& (\BH^{\bullet})^{\vee}\arrow[d] \arrow[r] & (\BE^{\bullet})^{\vee} \arrow[d] \arrow[r] &  \Omega_{B}[1] \arrow[d,equal] \\
& \BL_{M_Y} \arrow[r] &\BL_{M_{\CY/B}/B} \arrow[r] & \Omega_{B}[1]  
\end{tikzcd}
\]
By the same argument as in \cite[Section 2.3]{KT}, the composition 
\[(\BH^{\bullet})^{\vee}\rightarrow  (\BE^{\bullet})^{\vee} \rightarrow (\BE_{\mathrm{red}}^{\bullet})^{\vee}\]
is an isomorphism. This finishes the proof of the proposition.  
\qed
\\

For example, if $M=S^{[n]}$ and $\check{\beta}_1 \neq 0$, i.e.\ the curve class is not exceptional, we can use a commutative deformation given by the infinitesimal twistor family 
$\CS=\CY \rightarrow B$ with respect to the class $\check{\beta}_1$. The situation becomes more complicated already in the case of $S^{[n]}$ and $\check{\beta}_1 =0$ (i.e.\ an exceptional curve class). in this case, a commutative first-order deformation can no longer satisfy the property stated in Proposition \ref{nonc}. If $n=2$ and $S^{[2]}$ is isomorphic to a Fano variety of lines of some special cubic fourfold (e.g.\ see \cite[Theorem 1.0.3]{Ha}), then 
\[D_{\mathrm{perf}}(Y)=\langle D_{\mathrm{perf}}(S), \CO, \CO(1), \CO(2) \rangle,\]
and the family $\CY\rightarrow B$ is given by deformation of $Y$ away from the Hassett divisor.
\begin{remark} In \cite{To2}, Toda constructed geometric realisations of infinitesimal non-commutative deformations in $HH^2(X)$ for a smooth projective $X$. However, it is not clear, if they are of the type required by Proposition \ref{reduction}. In principle, there should be no problem in proving  Proposition \ref{reduction}, dropping the assumption. For that, one has to show that Toda's infinitesimal deformations behave well under base change.
\end{remark}

%
		

\bibliographystyle{amsalpha}
\bibliography{QMsfinal}

\providecommand{\bysame}{\leavevmode\hbox to3em{\hrulefill}\thinspace}
\providecommand{\MR}{\relax\ifhmode\unskip\space\fi MR }
\providecommand{\MRhref}[2]{%
  \href{http://www.ams.org/mathscinet-getitem?mr=#1}{#2}
}
\providecommand{\href}[2]{#2}
\begin{thebibliography}{MNOP06}

\bibitem[BF97]{BeFa}
K.~Behrend and B.~Fantechi, \emph{The intrinsic normal cone}, Invent. Math.
  \textbf{128} (1997), no.~1, 45--88.

\bibitem[BF03]{BF}
R.-O. Buchweitz and H.~Flenner, \emph{A semiregularity map for modules and
  applications to deformations}, Compositio Math. \textbf{137} (2003), no.~2,
  135--210.

\bibitem[Bri11]{Bri}
T.~Bridgeland, \emph{Hall algebras and curve-counting invariants}, J. Am. Math.
  Soc. \textbf{24} (2011), no.~4, 969--998.

\bibitem[C{\u{a}}l05]{Cal}
A.~C{\u{a}}ld{\u{a}}raru, \emph{The {M}ukai pairing. {II}. {T}he
  {H}ochschild-{K}ostant-{R}osenberg isomorphism}, Adv. Math. \textbf{194}
  (2005), no.~1, 34--66.

\bibitem[CKL17]{CKL}
H.-L. Chang, Y.-H. Kiem, and J.~Li, \emph{Torus localization and wall crossing
  for cosection localized virtual cycles}, Adv. Math. \textbf{308} (2017),
  964--986.

\bibitem[CKM14]{CFKM}
I.~{Ciocan-Fontanine}, B.~{Kim}, and D.~{Maulik}, \emph{{Stable quasimaps to
  GIT quotients}}, {J. Geom. Phys.} \textbf{75} (2014), 17--47.

\bibitem[Has00]{Ha}
B.~Hassett, \emph{Special cubic fourfolds}, Compositio Math. \textbf{120}
  (2000), no.~1, 1--23.

\bibitem[HT10]{HT}
D.~{Huybrechts} and R.~P. {Thomas}, \emph{{Deformation-obstruction theory for
  complexes via Atiyah and Kodaira-Spencer classes}}, {Math. Ann.} \textbf{346}
  (2010), no.~3, 545--569.

\bibitem[Huy99]{HuyK}
D.~Huybrechts, \emph{Compact hyperk{\"a}hler manifolds: {Basic} results},
  Invent. Math. \textbf{135} (1999), no.~1, 63--113.

\bibitem[JS12]{JS}
D.~Joyce and Y.~Song, \emph{A theory of generalized {Donaldson}-{Thomas}
  invariants}, Mem. Am. Math. Soc., vol. 1020, Providence, RI: American
  Mathematical Society (AMS), 2012.

\bibitem[KL13]{KiL}
Y.-H. Kiem and J.~Li, \emph{Localizing virtual cycles by cosections}, J. Amer.
  Math. Soc. \textbf{26} (2013), no.~4, 1025--1050.

\bibitem[KS08]{KS}
M.~Kontsevich and Y.~Soibelman, \emph{{Stability structures, motivic
  Donaldson-Thomas invariants and cluster transformations}}, arXiv:0811.2435
  (2008).

\bibitem[KT14]{KTred}
M.~Kool and R.~Thomas, \emph{Reduced classes and curve counting on surfaces.
  {I}: {Theory}}, Algebr. Geom. \textbf{1} (2014), no.~3, 334--383.

\bibitem[KT18]{KT}
M.~{Kool} and R.~P. {Thomas}, \emph{{Stable pairs with descendents on local
  surfaces. I: The vertical component}}, {Pure Appl. Math. Q.} \textbf{13}
  (2018), no.~4, 581--638.

\bibitem[LW15]{LW}
J.~Li and B.~Wu, \emph{Good degeneration of {Quot}-schemes and coherent
  systems}, Commun. Anal. Geom. \textbf{23} (2015), no.~4, 841--921.

\bibitem[MNOP06]{MNOP}
D.~Maulik, N.~Nekrasov, A.~Okounkov, and R.~Pandharipande,
  \emph{Gromov-{W}itten theory and {D}onaldson-{T}homas theory. {II}}, Compos.
  Math. \textbf{142} (2006), no.~5, 1286--1304. \MR{2264665}

\bibitem[MO09]{MO2}
D.~Maulik and A.~Oblomkov, \emph{Donaldson-{Thomas} theory of {{\(\mathcal
  A_n\times \mathbb{P}^{1}\)}}}, Compos. Math. \textbf{145} (2009), no.~5,
  1249--1276.

\bibitem[MPT10]{MPT}
D.~Maulik, R.~Pandharipande, and R.~P. Thomas, \emph{Curves on {{\(K3\)}}
  surfaces and modular forms}, J. Topol. \textbf{3} (2010), no.~4, 937--996.

\bibitem[Nes25]{N}
D.~Nesterov, \emph{Quasimaps to moduli spaces of sheaves}, Forum Math. Pi
  \textbf{13} (2025), e12.

\bibitem[Obe18a]{Ob18}
G.~Oberdieck, \emph{Gromov-{Witten} invariants of the {Hilbert} schemes of
  points of a {{\(K3\)}} surface}, Geom. Topol. \textbf{22} (2018), no.~1,
  323--437.

\bibitem[Obe18b]{OberG}
\bysame, \emph{On reduced stable pair invariants}, Math. Z. \textbf{289}
  (2018), no.~1-2, 323--353.

\bibitem[Obe19]{OP1}
\bysame, \emph{{Gromov-Witten theory of \(K3 \times \mathbb{P}^1\) and
  quasi-Jacobi forms}}, {Int. Math. Res. Not.} \textbf{2019} (2019), no.~16,
  4966--5011.

\bibitem[Obe21]{Ob21}
\bysame, \emph{{Marked relative invariants and GW/PT correspondences}},
  arXiv:2112.11949 (2021).

\bibitem[Obe22]{Ob22}
\bysame, \emph{{Holomorphic anomaly equations for the Hilbert scheme of points
  of a K3 surface}}, arXiv:2202.03361 (2022).

\bibitem[Obe24]{Ob}
\bysame, \emph{Multiple cover formulas for {{\(K3\)}} geometries,
  wall-crossing, and {Quot} schemes}, Geom. Topol. \textbf{28} (2024), no.~7,
  3221--3256.

\bibitem[OP10]{OP10b}
A.~Okounkov and R.~Pandharipande, \emph{Quantum cohomology of the {Hilbert}
  scheme of points in the plane}, Invent. Math. \textbf{179} (2010), no.~3,
  523--557.

\bibitem[OP16]{OPa}
G.~Oberdieck and R.~Pandharipande, \emph{Curve counting on {$K3\times E$}, the
  {I}gusa cusp form {$\chi_{10}$}, and descendent integration}, K3 surfaces and
  their moduli, Progr. Math., vol. 315, Birkh\"{a}user/Springer, 2016,
  pp.~245--278.

\bibitem[OP18]{OPi}
G.~Oberdieck and A.~{Pixton}, \emph{{Holomorphic anomaly equations and the
  Igusa cusp form conjecture}}, {Invent. Math.} \textbf{213} (2018), no.~2,
  507--587.

\bibitem[OS20]{OS2}
G.~Oberdieck and J.~{Shen}, \emph{{Curve counting on elliptic Calabi-Yau
  threefolds via derived categories}}, {J. Eur. Math. Soc.} \textbf{22} (2020),
  no.~3, 967--1002.

\bibitem[Per]{P}
A.~Perry, \emph{{The integral Hodge conjecture for two-dimensional Calabi-Yau
  categories}}, arXiv:2004.03163 (2020).

\bibitem[PT09]{PT}
R.~{Pandharipande} and R.~P. {Thomas}, \emph{{Curve counting via stable pairs
  in the derived category}}, {Invent. Math.} \textbf{178} (2009), no.~2,
  407--447.

\bibitem[Tod09]{To2}
Y.~Toda, \emph{Deformations and {F}ourier-{M}ukai transforms}, J. Differential
  Geom. \textbf{81} (2009), no.~1, 197--224.

\bibitem[{Zho}22]{YZ}
Y.~{Zhou}, \emph{{Quasimap wall-crossing for GIT quotients}}, {Invent. Math.}
  \textbf{227} (2022), no.~2, 581--660.

\end{thebibliography}

\end{document}